\documentclass[a4paper,12pt]{article}

\usepackage{latexsym}
\usepackage{amssymb}
\usepackage{theorem}
\usepackage{amsmath}
\usepackage{amscd}
\usepackage{graphicx}
\usepackage{url}
\usepackage[colorlinks=true]{hyperref}
\hypersetup{
setpagesize=false,
 bookmarksnumbered=true,%
 bookmarksopen=true,%
 colorlinks=true,%
 linkcolor=blue,
 citecolor=blue,
 urlcolor=blue,
}

\usepackage[usenames]{color}
\newtheorem{theorem}{Theorem}[section]
\newtheorem{corollary}[theorem]{Corollary}

\newtheorem{example}[theorem]{Example}
\newtheorem{proposition}[theorem]{Proposition}
\newtheorem{remark}[theorem]{Remark}
\newtheorem{definition}[theorem]{Definition}

\numberwithin{equation}{section}

\newcommand{\demo}{\par\noindent{\it Proof. \/}\ }
\newcommand{\enD}{\hfill $\Box$\vspace{3truemm} \par}
\newcommand{\R}{\mathbb{R}}

\newcommand{\ord}{\operatorname{ord}}
\newcommand{\bn}{\mbox{\boldmath $n$}}
\newcommand{\bt}{\mbox{\boldmath $t$}}
\newcommand{\bs}{\mbox{\boldmath $s$}}
\newcommand{\bb}{\mbox{\boldmath $b$}}
\newcommand{\ba}{\mbox{\boldmath $a$}}

\newcommand{\bx}{\mbox{\boldmath $x$}}

\newcommand{\bz}{\mbox{\boldmath $z$}}
\newcommand{\be}{\mbox{\boldmath $e$}}
\newcommand{\bmu}{\mbox{\boldmath $\mu$}}

\begin{document}

\title{Surfaces of revolution of frontals in the Euclidean space}

\author{Masatomo Takahashi and Keisuke Teramoto}

\date{\today}

\maketitle
\begin{abstract} 
For Legendre curves, we consider surfaces of revolution of frontals. 
The surface of revolution of a frontal can be considered as a framed base surface. 
We give the curvatures and basic invariants for surfaces of revolution by using the curvatures of Legendre curves.
Moreover, we give properties of surfaces of revolution with singularities and cones.
\end{abstract}

\renewcommand{\thefootnote}{\fnsymbol{footnote}}

\section{Introduction}

The surface of revolution is one of classical object in differential geometry (cf. \cite{Gray,Kenmotsu1,Kenmotsu2,Matins-Saji-Santos-Teramoto}). 
It has been known that if the profile curve (the plane curve) across the axis of revolution, then the surface of revolution has cone-type singularity.
On the other hand, the profile curve is not regular, then the surface of revolution must have singularities. 

In \cite{Kenmotsu1,Kenmotsu2}, the surfaces of revolution of regular curves are investigated. 
Especially, {Kenmotsu gave} concrete construction of the surface of revolution {with prescribed mean curvature}. 
In \cite{Matins-Saji-Santos-Teramoto}, the surfaces of revolution of singular curves are investigated. 
They construct a given the unbounded mean curvature of the surface of revolution. 

In this paper, we consider more general situation. 
We consider frontals (Legendre curves) as singular plane curves and framed base surfaces (framed surfaces) as singular surfaces. 
In \cite{Fukunaga-Takahashi1}, we give the curvature of Legendre curves in order to {analyze} Legendre curves. 
In \cite{Fukunaga-Takahashi3}, we give the basic invariants of framed surfaces {so as to analyze} framed surfaces. 
In \S 2, we review the theories of Legendre curves in the unit tangent bundle over the Euclidean plane {$\R^2$} 
and framed surfaces in the Euclidean space {$\R^3$}.
The surface of revolution of a frontal is a framed base surface.
We can deal with surfaces of revolution with singular points more directly.
In fact, we give the curvatures and basic invariants for surfaces of revolution by using the curvatures of Legendre curves in \S 3. 
Moreover, we give profile curves for given information of the curvatures, for instance, the Gauss curvature or the mean curvature. 
For the cases of constant Gauss and mean curvature surfaces of revolution, see \cite{Gray}.
We also consider relations between right-left equivalence relations among profile curves and right-left equivalence relations among the surface of revolution.
In particular, we give characterizations of $j/i$-cusps $((i,j)=(2,3),(2,5),(3,4),(3,5))$ appearing on profile curves of surfaces of revolution in terms of the  curvatures of Legendre curves. 
Further, applying such characterizations to the case of surfaces of revolution with constant Gauss or constant mean curvature. 
We give conditions that constant Gauss curvature surfaces of revolution 
have $j/i$-cusps as $((i,j)=(2,3),(3,4),(3,5))$ singularities by the construction data, 
and show that there are no constant Gauss curvature surfaces of revolution with $5/2$-cusp.
In contrast, we show that profile curves of constant mean curvature surfaces of revolution do not have $j/i$-cusps as singularities. 
Further, we classify cone-type singularities in some cases. 
In \S 4, we give concrete examples.

All maps and manifolds considered here are differentiable of class $C^{\infty}$ unless stated otherwise. 

\bigskip
\noindent
{\bf Acknowledgement}. 
The first author was partially supported by JSPS KAKENHI Grant Number JP 17K05238 and the second author was partially supported by  JSPS KAKENHI Grant Number JP 17J02151. 

\section{Preliminary}
We quickly review the theories of Legendre curves in the unit tangent bundle over the Euclidean plane (cf. \cite{Fukunaga-Takahashi1}) and framed surface in the Euclidean space (cf. \cite{Fukunaga-Takahashi3}).

\subsection{Legendre curves}

Let $\gamma:I\to\R^2$ and $\nu:I\to S^1$ be smooth mappings, where $I$ is an interval of $\R$ and $S^1$ is the unit circle.
We say that $(\gamma,\nu):I \to \R^2 \times S^1$ is {\it a Legendre curve} if $(\gamma,\nu)^*\theta=0$ for all $t \in I$, where $\theta$ is a canonical contact form on the unit tangent bundle $T_1 \R^2=\R^2 \times S^1$ over $\R^2$ (cf. \cite{Arnold1, Arnold2}). 
This condition is equivalent to $\dot{\gamma}(t) \cdot \nu(t)=0$ for all $t \in I$, 
where $\dot{\gamma}(t)=(d\gamma/dt)(t)$ and $\ba\cdot\bb=a_1b_1+a_2b_2$ for any $\ba=(a_1,a_2), \bb=(b_1,b_2)\in\R^2$. 
{A point $t_0\in I$ is called a {\it singular point} of $\gamma$ if $\dot{\gamma}(t_0)=0$.} 
When a Legendre curve $(\gamma,\nu):I\to\R^2\times S^1$ gives an immersion, it is called a {\it Legendre immersion}.
We say that $\gamma:I \to \R^2$ is a {\it frontal} (respectively, a {\it front}) if there exists $\nu:I \to S^1$ such that $(\gamma,\nu)$ is a Legendre curve (respectively, a Legendre immersion). 
Examples of Legendre curves see \cite{Ishikawa,Ishikawa-book}.
We have the Frenet type formula of a frontal $\gamma$ as follows.
We put on $\bmu (t)=J(\nu (t))$, where $J$ is the anticlockwise rotation by angle $\pi/2$ in $\R^2$.
Then $\{\nu(t), \bmu(t) \}$ is a moving frame of the frontal $\gamma(t)$ in $\R^2$ and we have the Frenet type formula,
\begin{equation}\label{Frenet.frontal}
\left(
\begin{array}{c}
\dot{\nu}(t)\\
\dot{\bmu}(t)
\end{array}
\right)
=
\left(
\begin{array}{cc}
0 & \ell(t)\\
-\ell(t) & 0
\end{array}
\right)
\left(
\begin{array}{c}
\nu(t)\\
\bmu(t)
\end{array}
\right), \ 
\dot\gamma(t) = \beta(t) \bmu(t),
\end{equation}
where $\ell(t)=\dot\nu(t) \cdot \bmu(t)$ and $\beta(t)=\dot{\gamma}(t) \cdot \bmu(t)$.
We call the pair $(\ell,\beta)$ {\it the curvature of the Legendre curve}. 
By \eqref{Frenet.frontal}, we see that $(\gamma,\nu)$ is a Legendre immersion if and only if $(\ell,\beta)\neq(0,0)$. 

\begin{definition}\label{congruent}{\rm
Let $(\gamma,\nu)$ and $(\widetilde{\gamma},\widetilde{\nu}):I \to \R^2 \times S^1$ be Legendre curves.
We say that $(\gamma,\nu)$ and $(\widetilde{\gamma},\widetilde{\nu})$ are {\it congruent as Legendre curves} if there exist a constant rotation $A \in SO(2)$ and a translation $\ba$ on $\R^2$ such that $\widetilde{\gamma}(t)=A(\gamma(t))+\ba$ and $\widetilde{\nu}(t)=A (\nu(t))$ for all $t \in I$. 
}
\end{definition}
\begin{theorem}[Existence Theorem for Legendre curves \cite{Fukunaga-Takahashi1}] \label{existence.Legendre}
Let $(\ell,\beta):I \to \R^2$ be a smooth mapping. 
There exists a Legendre curve $(\gamma,\nu):I \to \R^2 \times S^1$ whose curvature of the Legendre curve is $(\ell, \beta)$.
\end{theorem}
Actually, we have the following.
\begin{align*}
\gamma(t) &=\left(-\int \beta(t) \sin \left( \int \ell(t)\ dt\right) dt,\ \int \beta(t) \cos \left(\int \ell(t)\ dt\right) dt\right),\\
\nu(t) &= \left(\cos \left(\int \ell(t) \ dt\right), \ \sin \left(\int \ell(t)\ dt\right) \right).
\end{align*}

\begin{theorem}[Uniqueness Theorem for Legendre curves \cite{Fukunaga-Takahashi1}] \label{uniqueness.Legendre}
Let $(\gamma,\nu)$ and $(\widetilde{\gamma},\widetilde{\nu}):I \to \R^2 \times S^1$ be Legendre curves with the curvatures of Legendre curves $(\ell,\beta)$ and $(\widetilde{\ell},\widetilde{\beta})$, respectively.
Then $(\gamma,\nu)$ and $(\widetilde{\gamma},\widetilde{\nu})$ are congruent as Legendre curves if and only if $(\ell,\beta)$ and $(\widetilde{\ell},\widetilde{\beta})$ coincide.
\end{theorem}

Let $(\gamma,\nu):I \to \R^2 \times S^1$ be a Legendre curve with the curvature $(\ell,\beta)$. 
We define the {\it parallel curves} of $(\gamma,\nu)$ by $\gamma^\lambda:I \to \R^2, \gamma^\lambda(t)=\gamma(t)+\lambda \nu(t)$ for $\lambda \in \R$. 
Then $(\gamma^\lambda,\nu):I \to \R^2 \times S^1$ is also a Legendre curve with the curvature $(\ell,\beta+\lambda \ell)$.

Moreover, we define the {\it evolute} of $(\gamma,\nu)$ by $\mathcal{E}v(\gamma):I \to \R^2, 
\mathcal{E}v(\gamma)(t)=\gamma(t)+(\beta(t)/\ell(t))\nu(t)$, where $\ell(t) \not=0$ for all $t \in I$ (cf. \cite{Fukunaga-Takahashi2}).

\subsection{Framed surfaces}

Let $\R^3$ be the $3$-dimensional Euclidean space equipped with the inner product $\ba \cdot \bb = a_1 b_1 + a_2 b_2 + a_3 b_3$, 
where $\ba = (a_1, a_2, a_3)$ and $\bb = (b_1, b_2, b_3) \in \R^3$. 
The norm of $\ba$ is given by $\vert \ba \vert = \sqrt{\ba \cdot \ba}$. 
We also define the vector product
$$
\ba \times \bb=\det \left(
\begin{array}{ccc}
\be_1 & \be_2 & \be_3 \\
a_1 & a_2 & a_3 \\
b_1 & b_2 & b_3 
\end{array}
\right),
$$
where $\{ \be_1, \be_2, \be_3 \}$ is the canonical basis of $\R^3$. 
Let $U$ be a simply connected domain of $\R^2$ and $S^2$ be the unit sphere in $\R^3$, that is, $S^2=\{\ba \in \R^3| |\ba|=1\}$.
We denote a $3$-dimensional smooth manifold $\{(\ba,\bb) \in S^2 \times S^2| \ba \cdot \bb=0\}$ by $\Delta$.
\par
We say that $(\bx,\bn,\bs):U \to \R^3 \times \Delta$ is a {\it framed surface} if $\bx_u (u,v) \cdot \bn (u,v)=0$ and $\bx_v(u,v) \cdot \bn(u,v)=0$ for all $(u,v) \in U$, where $\bx_u(u,v)=(\partial \bx/\partial u)(u,v)$ and $\bx_v(u,v)=(\partial \bx/\partial v)(u,v)$. 
We say that $\bx:U \to \R^3$ is a {\it framed base surface} if there exists $(\bn,\bs):U \to \Delta$ such that $(\bx,\bn,\bs)$ is a framed surface.
\par
Similarly to the case of Legendre curves, 
the pair $(\bx,\bn):U\to\R^3\times S^2$ is said to be a {\it Legendre surface} 
if $\bx_u (u,v) \cdot \bn (u,v)=0$ and $\bx_v(u,v) \cdot \bn(u,v)=0$ for all $(u,v) \in U$. 
Moreover, when a Legendre surface $(\bx,\bn):U\to\R^3\times S^2$ gives an immersion, this is called a {\it Legendre immersion}. 
We say that $\bx:U\to\R^3$ be a {\it frontal} (respectively, a {\it front}) if there exists a map $\bn:U\to\ S^2$ such that 
the pair $(\bx,\bn):U\to\R^3\times S^2$ is a Legendre surface (respectively, a Legendre immersion).
By definition, the framed base surface is a frontal. 
At least locally, the frontal is a framed base surface. 
{For a framed surface $(\bx,\bn,\bs)$, we say that a point $p\in U$ is a {\it singular point of $\bx$} 
if $\bx$ is not an immersion at $p$.}

We denote $\bt(u,v)=\bn(u,v) \times \bs(u,v)$. 
Then $\{\bn(u,v),\bs(u,v),\bt(u,v)\}$ is a moving frame along $\bx(u,v)$.
Thus, we have the following systems of differential equations:
\begin{equation}\label{tangent}
\left(\begin{array}{c}
\bx_u \\
\bx_v
\end{array}\right)
=
\left(\begin{array}{cc} 
a_1 & b_1 \\
a_2 & b_2 
\end{array}\right)
\left(\begin{array}{c}
\bs \\
\bt
\end{array}\right),
\end{equation}
\begin{equation}\label{frame}
\left(\begin{array}{c} 
\bn_u \\
\bs_u \\
\bt_u
\end{array}\right)
=
\left(\begin{array}{ccc} 
0 & e_1 & f_1 \\
-e_1 & 0 & g_1 \\
-f_1 & -g_1 & 0
\end{array}\right)
\left(\begin{array}{c} 
\bn \\
\bs \\
\bt
\end{array}\right)
, \ 
\left(\begin{array}{c} 
\bn_v \\
\bs_v \\
\bt_v
\end{array}\right)
=
\left(\begin{array}{ccc} 
0 & e_2 & f_2 \\
-e_2 & 0 & g_2 \\
-f_2 & -g_2 & 0
\end{array}\right)
\left(\begin{array}{c} 
\bn \\
\bs \\
\bt
\end{array}\right),
\end{equation}
where $a_i,b_i,e_i,f_i,g_i:U \to \R, i=1,2$ are smooth functions and we call the functions {\it basic invariants} of the framed surface. 
We denote the matrices $(\ref{tangent})$ and $(\ref{frame})$ by $\mathcal{G}, \mathcal{F}_1, \mathcal{F}_2$, respectively. 
We also call the matrices $(\mathcal{G}, \mathcal{F}_1, \mathcal{F}_2)$  {\it basic invariants} of the framed surface $(\bx,\bn,\bs)$.
Since the integrability condition $\bx_{uv}=\bx_{vu}$ and $\mathcal{F}_{2,u}-\mathcal{F}_{1,v}
=\mathcal{F}_1\mathcal{F}_2-\mathcal{F}_2\mathcal{F}_1$, 
the basic invariants should satisfy the following conditions:
\begin{equation}\label{integrability.condition}
\begin{cases}
a_{1,v}-b_1g_2 = a_{2,u}-b_2g_1, \\
b_{1,v}-a_2g_1 = b_{2,u}-a_1g_2, \\
a_1 e_2 + b_1 f_2 = a_2 e_1 + b_2 f_1,
\end{cases}
\begin{cases}
e_{1,v}-f_1g_2 = e_{2,u}-f_2g_1, \\
f_{1,v}-e_2g_1 = f_{2,u}-e_1g_2, \\
g_{1,v}-e_1f_2 = g_{2,u}-e_2f_1. 
\end{cases}
\end{equation}

We give fundamental theorems for framed surfaces, that is, 
the existence and uniqueness theorem of framed surfaces for basic invariants.

\begin{definition}\label{congruent.framed.surface}{\rm
Let $(\bx,\bn,\bs), (\widetilde{\bx},\widetilde{\bn},\widetilde{\bs}):U \to \R^3 \times \Delta$ be framed surfaces.
We say that $(\bx,\bn,\bs)$ and $(\widetilde{\bx},\widetilde{\bn},\widetilde{\bs})$ are {\it congruent as framed surfaces} if there exist a constant rotation $A \in SO(3)$ and a translation $\ba \in \R^3$ such that 
$$
\widetilde{\bx}(u,v)=A(\bx(u,v))+\ba, \ \widetilde{\bn}(u,v)=A(\bn(u,v)), \  \widetilde{\bs}(u,v)=A(\bs(u,v))
$$ 
for all $(u,v) \in U$.}
\end{definition}

\begin{theorem}[Existence Theorem for framed surfaces \cite{Fukunaga-Takahashi3}]\label{existence.framed.surface}
Let $U$ be a simply connected domain in $\R^2$ and let $a_i,b_i,e_i,f_i,g_i:U \to \R, i=1,2$ be smooth functions with the integrability conditions $(\ref{integrability.condition})$.
Then there exists a framed surface $(\bx,\bn,\bs):U \to \R^3 \times \Delta$ whose associated basic invariants is $a_i,b_i,e_i,f_i,g_i, i=1,2$.
\end{theorem}
\begin{theorem}[Uniqueness Theorem for framed surfaces \cite{Fukunaga-Takahashi3}]\label{uniqueness.framed.surface}
Let $(\bx,\bn,\bs)$, $(\widetilde{\bx},\widetilde{\bn},\widetilde{\bs}):U \to \R^3 \times \Delta$ be framed surfaces 
with basic invariants $(\mathcal{G},\mathcal{F}_1,\mathcal{F}_2)$ and $(\widetilde{\mathcal{G}},\widetilde{\mathcal{F}}_1,\widetilde{\mathcal{F}}_2)$, respectively.
Then $(\bx,\bn,\bs)$ and $(\widetilde{\bx},\widetilde{\bn},\widetilde{\bs})$ are congruent as framed surfaces if and only if $(\mathcal{G},\mathcal{F}_1,\mathcal{F}_2)$ and $(\widetilde{\mathcal{G}},\widetilde{\mathcal{F}}_1,\widetilde{\mathcal{F}}_2)$  coincide.
\end{theorem}

Let $(\bx,\bn,\bs):U \to \R^3 \times \Delta$ be a framed surface with basic invariants $(\mathcal{G},\mathcal{F}_1,\mathcal{F}_2).$

\begin{definition}\label{curvature.framed.surface}{\rm 
We define a smooth mapping $C_F=(J_F,K_F,H_F):U \to \R^3$ by 
\begin{align*}
J_F =
\det \left(\begin{array}{cc} 
a_1 & b_1 \\ 
a_2 & b_2
\end{array}\right), 
K_F =
\det \left(\begin{array}{cc} 
e_1 & f_1 \\ 
e_2 & f_2
\end{array}\right), 
H_F = -\frac{1}{2}\left\{ 
\det \left(\begin{array}{cc} 
a_1 & f_1 \\ 
a_2 & f_2
\end{array}\right)
-
\det \left(\begin{array}{cc} 
b_1 & e_1 \\ 
b_2 & e_2
\end{array}\right)
\right\}.
\end{align*}
We call $C_F = (J_F,K_F,H_F)$ a {\it curvature of the framed surface}.
\par
We also define a smooth mapping $I_F:U \to \R^8$ by 
\begin{multline*}
I_F=\left(C_F,
\det \left(\begin{array}{cc} 
a_1 & g_1 \\ 
a_2 & g_2
\end{array}\right),
\det \left(\begin{array}{cc} 
b_1 & g_1 \\ 
b_2 & g_2
\end{array}\right),\right.\\ \left.
\det \left(\begin{array}{cc} 
e_1 & g_1 \\ 
e_2 & g_2
\end{array}\right), 
\det \left(\begin{array}{cc} 
f_1 & g_1 \\ 
f_2 & g_2
\end{array}\right),
\det \left(\begin{array}{cc} 
a_1 & e_1 \\ 
a_2 & e_2
\end{array}\right)
\right).
\end{multline*}
We call the mapping $I_F:U \to \R^8$ a {\it concomitant mapping} of the framed surface $(\bx,\bn,\bs)$. 
}
\end{definition}
\begin{remark}{\rm 
If the surface $\bx$ is regular, then we have $K=K_F/J_F$ and $H=H_F/J_F$, where 
$K$ is the Gauss curvature and $H$ is the mean curvature of the regular surface (cf. \cite{Fukunaga-Takahashi3}). 
For relations between behaviour of the Gauss curvature, the mean curvature of fronts at non-degenerate singular points and geometric invariants of fronts see \cite{Martins-Saji-Umehara-Yamada}.} 
\end{remark}
We say that $(\bx,\bn,\bs):U \to \R^3 \times \Delta$ is a {\it framed immersion} if $(\bx,\bn,\bs)$ is an immersion.
\begin{proposition}[\cite{Fukunaga-Takahashi3}]
\label{immersive.condition}
Let $(\bx,\bn,\bs):U \to \R^3 \times \Delta$ be a framed surface and $p \in U$.
\par
$(1)$ $\bx$ is an immersion (a regular surface) around $p$ if and only if $J_F (p) \not=0$.
\par
$(2)$ $(\bx,\bn)$ is a Legendre immersion around $p$ if and only if $C_F (p) \not=0$.
\par
$(3)$ $(\bx,\bn,\bs)$ is a framed immersion around $p$ if and only if $I_F (p) \not=0$.
\end{proposition}

Let $(\bx,\bn,\bs):U \to \R^3 \times \Delta$ be a framed surface with basic invariants $(\mathcal{G},\mathcal{F}_1,\mathcal{F}_2)$. 
We define the {\it parallel surface} $\bx^\lambda : U \rightarrow \mathbb{R}^3$ of the framed surface $(\bx,\bn,\bs)$ by  $\bx^\lambda(u,v) = \bx(u,v) + \lambda \bn(u,v)$, where $\lambda \in \mathbb{R}$. 
Then $(\bx^\lambda,\bn,\bs):U \to \R^3 \times \Delta$ is also a framed surface and  basic invariants $(\mathcal{G}^\lambda,  \mathcal{F}_1^\lambda, \mathcal{F}_2^\lambda)$ are given by  
\begin{equation*}
\mathcal{G}^\lambda = \mathcal{G} + \lambda \left(\begin{array}{cc} e_1 & f_1 \\ e_2 & f_2\end{array}\right), \ \mathcal{F}_1^\lambda =
\mathcal{F}_1, \ \mathcal{F}_2^\lambda =\mathcal{F}_2.
\end{equation*}
By a direct calculation, we have the curvature $(J^\lambda_F,K^\lambda_F,H^\lambda_F)$ of the framed surface $(\bx^\lambda,\bn,\bs)$ is given by 
$$
J^\lambda_F=J_F-2H_F\lambda+K_F\lambda^2, \ K^\lambda_F=K_F, \  H^\lambda_F=H_F-K_F\lambda.
$$

We also define the {\it evolute} (or, the {\it focal surface}) of the framed surface $(\bx,\bn,\bs)$ by $\mathcal{E}v(\bx):U \to \R^3, \mathcal{E}v(\bx)(u,v)=\bx(u,v)+\lambda \bn(u,v)$, where 
$\lambda$ is a solution of the equation
$$
K_F(u,v) \lambda^2-2H_F(u,v) \lambda+J_F(u,v)=0
$$
for all $(u,v) \in U$. 
The study of focal surfaces of fronts from a different viewpoint is known in \cite{Teramoto2}.

\begin{remark}\label{rem:evolute1}{\rm
Even if $K_F(u,v)=0$ for all $(u,v) \in U$, we can define an evolute when $H_F(u,v) \not=0$. 
In this case, we have only one evolute. 
}
\end{remark}

Moreover, if we consider a similar surface $(r\bx,\bn,\bs):U \to \R^3 \times \Delta$ for non-zero constant $r \in \R$, 
Then the curvature $(J^r_F,K^r_F,H^r_F)$ of the framed surface $(r\bx,\bn,\bs)$ is given by 
$$
J^r_F=r^2 J_F, \ K^r_F=K_F, \  H^r_F=r H_F.
$$

\section{Surfaces of revolution of frontals}

Let $(\gamma,\nu):I \to \R^2 \times S^1$ be  a Legendre curve with the curvature $(\ell,\beta)$. 
We denote $\gamma(t)=(x(t),z(t))$ and $\nu(t)=(a(t),b(t))$. 
By definition and the Frenet type formula (\ref{Frenet.frontal}), we have $\dot{x}(t)a(t)+\dot{z}(t)b(t)=0$, $a^2(t)+b^2(t)=1$, 
\begin{equation}\label{Frenet-type}
\left(\begin{array}{c} 
\dot{x}(t) \\
\dot{z}(t)
\end{array}\right)=
\beta(t)
\left(\begin{array}{c} 
-b(t) \\
a(t)
\end{array}\right), 
\left(\begin{array}{c} 
\dot{a}(t) \\
\dot{b}(t)
\end{array}\right)=
\ell(t)
\left(\begin{array}{c} 
-b(t) \\
a(t)
\end{array}\right)
\end{equation}
for all $t \in I$. 
Set a smooth function $\varphi:I \to \R$ such that $a(t)=\cos \varphi(t)$ and $b(t)=\sin \varphi(t)$. 

We consider $\R^2 \subset \R^3$ as $(x,z)$-plane into $(x,y,z)$-space. 
We give the two surfaces of revolution, that is, around $x$-axis and $z$-axis of the frontal $\gamma$, respectively. 
We call $\gamma$ {\it a profile curve} (cf. \cite{Gray}).

First we consider the surface of revolution around $x$-axis.
We denote the {\it surface of revolution of $\gamma$ around $x$-axis} by $\bx:I \times [0,2\pi) \to \R^3$,
$$
\bx(t,\theta)=(x(t),z(t) \cos \theta, z(t) \sin \theta).
$$
Then the surface of revolution around $x$-axis of the frontal is a framed base surface. 
\begin{proposition}\label{x-axis.rotation}
Under the above notation, 
$(\bx,\bn^x,\bs^x):I \times [0,2\pi) \to \R^3 \times \Delta$ is a framed surface with basic invariants,
\begin{equation}\label{x-axis.rotation.invariant}
\mathcal{G}^x=
\left(\begin{array}{cc} 
0 & -\beta(t)  \\
-z(t) & 0
\end{array}\right),
\mathcal{F}^x_1=
\left(\begin{array}{ccc} 
0 & 0& \ell(t)\\
0 & 0& 0\\
-\ell(t) & 0& 0
\end{array}\right),
\mathcal{F}^x_2=
\left(\begin{array}{ccc} 
0 & b(t)& 0\\
-b(t) & 0& -a(t)\\
0 & a(t)& 0
\end{array}\right),
\end{equation}
where 
$\bn^x(t,\theta)=(-a(t),-b(t) \cos \theta,-b(t) \sin \theta), 
\bs^x(t,\theta)=(0,\sin \theta, -\cos \theta)$.
\end{proposition}
\demo
By a direct calculation, 
$\bt^x(t,\theta)=\bn^x(t,\theta) \times \bs^x(t,\theta)=(b(t),-a(t)\cos  \theta, -a(t)\sin \theta)$.
Since 
\begin{align*}
\bx_t(t,\theta)&= (-\beta(t)b(t),\beta(t)a(t)\cos \theta,\beta(t)a(t)\sin \theta)=-\beta(t)\bt^x(t,\theta),\\ 
\bx_\theta(t,\theta)&=(0,-z(t)\sin \theta,z(t)\cos \theta)=-z(t)\bs^x(t,\theta),\\
\bn^x_t(t,\theta)&=(\ell(t)b(t),-\ell(t)a(t)\cos \theta,-\ell(t)a(t)\sin \theta)=\ell(t)\bt^x(t,\theta),\\
\bn^x_{\theta}(t,\theta)&=(0,b(t)\sin \theta,-b(t)\cos \theta)=b(t)\bs^x(t,\theta),\\
\bs^x_t(t,\theta)&=0,\\
\bs^x_{\theta}(t,\theta)&=(0,\cos \theta,\sin \theta)=-b(t)\bn^x(t,\theta)-a(t)\bt^x(t,\theta),\\
\bt^x_{t}(t,\theta)&=(\ell(t)a(t), \ell(t)b(t)\cos \theta, \ell(t)b(t)\sin \theta)=-\ell(t)\bn^x(t,\theta),\\
\bt^x_{\theta}(t,\theta)&=(0,b(t)\sin \theta,b(t)\cos \theta)=a(t)\bs^x(t,\theta),
\end{align*}
we have basic invariants $(\mathcal{G}^x,\mathcal{F}^x_1,\mathcal{F}^x_2)$.
\enD
By a direct calculation, we have 
\begin{align*}
&J^x_F(t,\theta)=-\beta(t)z(t), K^x_F(t,\theta)=-b(t)\ell(t), H^x_F(t,\theta)=-\frac{1}{2}(z(t)\ell(t)+\beta(t)b(t)),\\
&\det \left(\begin{array}{cc} 
a_1 & g_1 \\ 
a_2 & g_2
\end{array}\right)(t,\theta)=0,
\det \left(\begin{array}{cc} 
b_1 & g_1 \\ 
b_2 & g_2
\end{array}\right)(t,\theta)=\beta(t)a(t),
\det \left(\begin{array}{cc} 
e_1 & g_1 \\ 
e_2 & g_2
\end{array}\right)(t,\theta)=0,\\
&\det \left(\begin{array}{cc} 
f_1 & g_1 \\ 
f_2 & g_2
\end{array}\right)(t,\theta)=-\ell(t)a(t),
\det \left(\begin{array}{cc} 
a_1 & e_1 \\ 
a_2 & e_2
\end{array}\right)(t,\theta)=0.
\end{align*}

Next, we consider the surface of revolution around  $z$-axis.
We denote the {\it surface of revolution of $\gamma$ around $z$-axis} by $\bz:I \times [0,2\pi) \to \R^3$,
$$
\bz(t,\theta)=(x(t) \cos \theta,x(t) \sin \theta, z(t)).
$$
Then the surface of revolution around  $z$-axis of the frontal is also a framed base surface. 
By the similar calculation in Proposition \ref{x-axis.rotation}, we have the following result.
\begin{proposition}\label{z-axis.rotation}
Under the above notation, 
$(\bz,\bn^z,\bs^z):I \times [0,2\pi) \to \R^3 \times \Delta$ is a framed surface with  basic invariants,
\begin{equation}\label{z-axis.rotation.invariant}
\mathcal{G}^z=
\left(\begin{array}{cc}
0 & -\beta(t)  \\
-x(t) & 0
\end{array}\right),
\mathcal{F}^z_1=
\left(\begin{array}{ccc}
0 & 0& -\ell(t)\\
0 & 0& 0\\
\ell(t) & 0& 0
\end{array}\right),
\mathcal{F}^z_2=
\left(\begin{array}{ccc}
0 & -a(t)& 0\\
a(t) & 0& b(t)\\
0 & -b(t)& 0
\end{array}\right),
\end{equation}
where 
$\bn^z(t,\theta)=(a(t) \cos \theta,a(t)\sin \theta,b(t)), 
\bs^z(t,\theta)=(\sin \theta,-\cos \theta,0)$.
\end{proposition}

By a direct calculation, we have 
\begin{align*}
&J^z_F(t,\theta)=-\beta(t)x(t), K^z_F(t,\theta)=-a(t)\ell(t), H^z_F(t,\theta)=\frac{1}{2}(x(t)\ell(t)+\beta(t)a(t)),\\
&\det \left(\begin{array}{cc} 
a_1 & g_1 \\ 
a_2 & g_2
\end{array}\right)(t,\theta)=0,
\det \left(\begin{array}{cc} 
b_1 & g_1 \\ 
b_2 & g_2
\end{array}\right)(t,\theta)=-\beta(t)b(t),
\det \left(\begin{array}{cc} 
e_1 & g_1 \\ 
e_2 & g_2
\end{array}\right)(t,\theta)=0,\\
&\det \left(\begin{array}{cc} 
f_1 & g_1 \\ 
f_2 & g_2
\end{array}\right)(t,\theta)=-\ell(t)b(t),
\det \left(\begin{array}{cc} 
a_1 & e_1 \\ 
a_2 & e_2
\end{array}\right)(t,\theta)=0.
\end{align*}
\begin{proposition} \label{front-frontal}
Under the above notation, we have the following.
\par
$(1)$ The surfaces of revolution $\bx$ and $\bz$ are frontals if and only if $\gamma$ is a frontal. 
\par
$(2)$ The surface of revolution $\bx$ or $\bz$ is a front if and only if $\gamma$ is a front. 
\end{proposition}
\demo
$(1)$ If $\bx$ and $\bz$ are frontals, then there exists a point $(t,\theta)$ such that $C^x_F(t,\theta)=0$ and $C^z_F(t,\theta)=0$ by Proposition \ref{immersive.condition}. 
Since $K^x_F(t,\theta)=K^z_F(t,\theta)=0$, we have $\ell(t)=0$.
Moreover, since $H^x_F(t,\theta)=H^z_F(t,\theta)=0$, we have $\beta(t)=0$. 
It follows that $\gamma$ is a frontal. 
\par
Conversely, let $\gamma$ be a frontal.
If $t$ is a singular point of $(\gamma,\nu)$, then $\ell(t)=\beta(t)=0$. 
By a direct calculation, we have $C^x_F(t,\theta)=0$ and $C^z_F(t,\theta)=0$. 
Therefore, $\bx$ and $\bz$ are frontals.
\par
$(2)$ By (1), we have the result.
\enD

\begin{proposition} \label{revolution.congruent}
$(\bx,\bn^x,\bs^x)$ and $(\bz,\bn^z,\bs^z)$ are congruent as framed surfaces if and only if the Legendre curve $(\gamma,\nu):I \to \R^2 \times S^1$ is given by 
$$
\gamma(t)=\left( \pm \int \frac{1}{\sqrt{2}} \beta(t) dt, \pm \int \frac{1}{\sqrt{2}} \beta(t) dt\right), 
\nu(t)=\left(\pm \frac{1}{\sqrt{2}}, \mp \frac{1}{\sqrt{2}}\right).
$$
\end{proposition}
\demo
By Theorem \ref{uniqueness.framed.surface}, 
$(\bx,\bn^x,\bs^x)$ and $(\bz,\bn^z,\bs^z)$ are congruent as framed surface if and only if $(\mathcal{G}^x,\mathcal{F}^x_1,\mathcal{F}^x_2)(t,\theta)=(\mathcal{G}^z,\mathcal{F}^z_1,\mathcal{F}_2)(t,\theta)$ for all $(t,\theta) \in I \times[0,2\pi)$.
It follows that $z(t)=x(t), \ell(t)=0, a(t)=-b(t)$ for all $t \in I$.
Since $\ell(t)=0$, $a(t)$ and $b(t)$ are constants.
By $a^2(t)+b^2(t)=1$, we have $\nu(t)=(\pm {1}{/\sqrt{2}}, \mp 1/\sqrt{2})$.
Moreover, $\dot{x}(t)=-\beta(t) b(t), \dot{z}(t)=\beta(t) a(t)$, we have $\gamma(t)=(x(t),z(t))$ by integration.
\enD


\begin{theorem}\label{existence.surface-revolution}
Suppose that $(\mathcal{G}^x,\mathcal{F}^x_1,\mathcal{F}^x_2)$ and $(\mathcal{G}^z,\mathcal{F}^z_1,\mathcal{F}^z_2)$ are given by the forms of $(\ref{x-axis.rotation.invariant})$ and of $(\ref{z-axis.rotation.invariant})$ with 
$\dot{x}(t)=-\beta(t)b(t), \dot{y}(t)=\beta(t)a(t), \dot{a}(t)=-\ell(t)b(t), \dot{b}(t)=\ell(t)a(t), a^2(t)+b^2(t)=1$. 
Then there exists a unique Legendre curve $(\gamma,\nu):I \to \R^2 \times S^1$ with the curvature $(\ell,\beta)$ such that basic invariants of surfaces of revolution around $x$-axis and $z$-axis are $(\mathcal{G}^x,\mathcal{F}^x_1,\mathcal{F}^x_2)$ and $(\mathcal{G}^z,\mathcal{F}^z_1,\mathcal{F}^z_2)$, respectively.
\end{theorem}
\demo
By the forms of (\ref{x-axis.rotation.invariant}) and of (\ref{z-axis.rotation.invariant}), we define $(\gamma,\nu):I \to \R^2 \times S^1$  by $\gamma(t)=(x(t),z(t))$ and $\nu(t)=(a(t),b(t))$.
Then $(\gamma,\nu)$ is a Legendre curve with the curvature $(\ell,\beta)$. 
Moreover, basic invariants of the surfaces of $x$-axis revolution $(\bx,\bn^x,\bs^x)$ and $z$-axis revolution $(\bz,\bn^z,\bs^z)$ of the Legendre curve $(\gamma,\nu)$ are $(\mathcal{G}^x,\mathcal{F}^x_1,\mathcal{F}^x_2)$ and $(\mathcal{G}^z,\mathcal{F}^z_1,\mathcal{F}^z_2)$, respectively. 
\enD

We now concentrate on the case of surfaces of revolution around $z$-axis. 
Let $(\gamma,\nu):I \to \R^2 \times S^1$ be a Legendre curve with the curvature $(\ell,\beta)$. We denote $\gamma(t)=(x(t),z(t))$ and $\nu(t)=(a(t),b(t))$. 

The mapping $f:U \to \R^3$ is a {\it cone} if $f$ is given by $f(u,v)=(u \cos v,u \sin v,u)$.
\begin{proposition} \label{cone}
If $x(t_0)=0, \beta(t_0) \not= 0, a(t_0) \not=0$ and $b(t_0) \not=0$, then the surface of revolution $\bz:I \times [0,2\pi) \to \R^3$ at $(t_0,\theta_0) \in I \times [0,2\pi)$ is diffeomorphic to the cone at $(0,\theta_0)$.
\end{proposition}
\demo
By the assumption, we have $\dot{x}(t_0) \not=0$ and $\dot{z}(t_0) \not=0$. 
The surface of revolution $\bz:I \times [0,2\pi) \to \R^3$ is given by 
$\bz(t,\theta)=(x(t) \cos \theta,x(t) \sin \theta, z(t)).$
By using the diffeomorphisms of the source $(t,\theta) \mapsto (x(t),\theta)$ and the target $(X,Y,Z) \mapsto (X,Y,z(Z))$, then 
$\bz(t,\theta)$ at $(t_0,\theta_0)$ is diffeomorphic to the cone $(t,\theta) \mapsto (t\cos \theta,t\sin \theta, t)$ at $(0,\theta_0)$.
\enD
For more general cases see Theorem \ref{diffeo-zero}.

\begin{proposition} \label{zero-Gauss-curvature}
If $K^z_F(t,\theta)=0$ for all $(t,\theta) \in I \times [0,2\pi)$, then $\gamma$ is a part of a line. 
\end{proposition}
\demo
If $K^z_F(t,\theta)=0$ for all $(t,\theta) \in I \times [0,2\pi)$, then $-a(t)\ell(t)=0$ for all $t \in I$.
By differentiating, we have $b(t)\ell^2(t)+a(t)\dot{\ell}(t)=0$.
If $\ell(t_0) \not =0$ at a point $t_0 \in I$, then $a(t)=0$ and $b(t)=0$ around $t_0$.
This is a contradict the fact that $(a(t),b(t)) \in S^1$. 
It follows that $\ell(t)=0$ for all $t \in I$. 
This is equivalent to $a(t)$ and $b(t)$ are constants. 
Hence, $\gamma$ is a part of a line.
\enD
By Proposition \ref{zero-Gauss-curvature}, the surfaces of revolution of a frontal {satisfying} 
$K^z_F(t,\theta)=0$ for all $(t,\theta) \in I \times [0,2\pi)$ are given by a part of a cone, a cylinder, a plane, a line, a circle or a point. 

We consider general cases. 
We denote $J(t)=J^z_F(t,\theta), K(t)=K^z_F(t,\theta)$ and $H(t)=H^z_F(t,\theta)$.
\begin{theorem}\label{Gauss-curvature}
{Let $(\gamma,\nu):(I,t_0) \to \R^2 \times S^1$ be a Legendre curve of the form $\gamma(t)=(x(t),z(t))$ and $\nu(t)=(\cos \varphi(t), \sin \varphi(t))$} 
with the curvature $(\dot{\varphi},\beta)$. 
\par
$(1)$ Suppose that $\beta$ is a real analytic around $t_0$ and there exists a  function $\alpha:(I,t_0) \to \R$ such that $K(t)=\alpha(t)J(t)$ and $\alpha(t)\beta^2(t)(t-t_0)^2$ is a real analytic around $t_0$. 
Then $x(t)$ is a solution of 
\begin{equation}\label{differential-eq}
\beta(t)\ddot{x}(t)-\dot{\beta}(t)\dot{x}(t)+\alpha(t)\beta^3(t)x(t)=0
\end{equation}
around $t_0$, 
\begin{equation*}
z(t) = \pm \int \beta(t) \left(1-\left(\int \alpha(t) \beta(t) x(t) dt\right)^2 \right)^\frac{1}{2} dt
\end{equation*}
and 
$$
\cos \varphi(t)= \pm \left(1-\left(\int \alpha(t)\beta(t)x(t) dt\right)^2\right)^{\frac{1}{2}}, \ \sin \varphi(t)=\int \alpha(t)\beta(t)x(t) dt.
$$
\par
$(2)$ Suppose that $x(t_0)>0$ and given smooth functions $J$ and $K: (I,t_0) \to \R$. 
Then  $(\gamma,\nu)$ is given by 
\begin{equation*}
x(t) = \left(-2 \int J(t) \left(\int K(t) dt\right) dt \right)^\frac{1}{2}, \ 
z(t) = \mp \int \frac{J(t)}{x(t)} \left(1-\left(\int K(t) dt\right)^2 \right)^\frac{1}{2} dt
\end{equation*}
and 
$$
\cos \varphi(t)=\pm \left(1-\left(\int K(t) dt\right)^2\right)^\frac{1}{2}, \  \sin \varphi(t)=-\int K(t) dt.
$$
\end{theorem}
\demo
$(1)$ Since $K(t)=\alpha(t)J(t)$, we have $\dot{\varphi}(t)\cos \varphi(t)=\alpha(t)\beta(t)x(t)$.
{By 
(\ref{Frenet-type})}, $\dot{x}(t)=-\beta(t)\sin \varphi(t)$. 
It follows that $\ddot{x}(t)=-\dot{\beta}(t)\sin \varphi(t)-\alpha(t)\beta^2(t)x(t)$ and $\beta(t) \ddot{x}(t)=\dot{\beta}(t)\dot{x}(t)-\alpha(t)\beta^3(t)x(t)$. 
Then $x(t)$ is satisfied a second order ordinary linear differential equation of $(\ref{differential-eq})$. 
By the assumption, $t_0$ is a regular singularity of $(\ref{differential-eq})$. 
Then there exists a solution of $x(t)$ around $t_0$ by using the method of Frobenius (cf. \cite{Tenenbaum-Pollard}). 
Moreover, $(d/dt)(\sin \varphi(t))=\dot{\varphi}(t)\cos \varphi(t)=\alpha(t)\beta(t)x(t)$, we have 
$$
\cos \varphi(t)= \pm \left(1-\left(\int \alpha(t)\beta(t)x(t) dt\right)^2\right)^{\frac{1}{2}},\ \sin \varphi(t)=\int \alpha(t)\beta(t)x(t) dt.
$$
Since $\dot{z}(t)=\beta(t) \cos \varphi(t)$, we have 
$$
z(t)=\int \beta(t)\cos \varphi(t) dt =\pm \int \beta(t) \left(1-\left(\int \alpha(t) \beta(t) x(t) dt\right)^2 \right)^\frac{1}{2} dt.
$$
\par
$(2)$ We have $J(t)=-\beta(t)x(t)$ and $K(t)=-\dot{\varphi}(t)\cos \varphi(t)$. 
By $(d/dt)(\sin \varphi(t))=\dot{\varphi}(t)\cos \varphi(t)=-K(t)$, we have 
$$
\cos \varphi(t)=\pm \left(1-\left(\int K(t) dt\right)^2\right)^\frac{1}{2}, \ \sin \varphi(t)=-\int K(t) dt.
$$
By the Frenet type formula (\ref{Frenet-type}), $\dot{x}(t)=-\beta(t)\sin \varphi(t)$. 
It follows that 
$$
\frac{d}{dt} x^2(t)=2 x(t) \dot{x}(t)=-2\beta(t)x(t)\sin\varphi(t)
=-2 J(t) \int K(t) dt.
$$
Then 
$$
x^2(t)=-2 \int J(t) \left(\int K(t) dt\right) dt.
$$
Since $x(t_0)>0$, we have $x(t)$ around $t_0$. 
{Further}, by $\dot{z}(t)=\beta(t)\cos \varphi(t)=-(J(t)/x(t)) \cos \varphi(t)$, we have
$$
z(t)=\mp \int \frac{J(t)}{x(t)} \left(1-\left(\int K(t) dt\right)^2\right)^\frac{1}{2} dt.
$$
\enD

By the similar calculation in \cite{Kenmotsu1,Matins-Saji-Santos-Teramoto}, we have the following result for $H$. 
\begin{theorem}\label{mean-curvature}
{Let $(\gamma,\nu):(I,t_0) \to \R^2 \times S^1$ be a Legendre curve of the form $\gamma(t)=(x(t),z(t))$ and $\nu(t)=(\cos \varphi(t), \sin \varphi(t))$} 
with the curvature $(\dot{\varphi},\beta)$. 
Suppose that $x(t_0)>0$, we give $\beta$ and there exists a smooth function $\alpha:(I,t_0) \to \R$ such that $H(t)=\alpha(t)J(t)$.
Then $(\gamma,\nu)$ is given by
\begin{equation}\label{mean-xz}
x(t) = (F(t)^2+G(t)^2)^\frac{1}{2}, \ 
z(t) = \int \frac{\beta(t)}{x(t)} \left( F(t)\sin \eta(t)-G(t)\cos \eta(t)\right) dt
\end{equation}
and 
$$
\cos \varphi(t)=\frac{F(t)\sin \eta(t)-G(t)\cos \eta(t)}{x(t)},\ \sin \varphi(t)=\frac{F(t)\cos \eta(t)+G(t)\sin \eta(t)}{x(t)},
$$
where 
$$
F(t)=-\int \beta(t)\cos \eta(t) dt, \ G(t)=-\int \beta(t)\sin \eta(t) dt, \  \eta(t)=2 \int \alpha(t)\beta(t) dt.
$$
\end{theorem}
\demo
Since $H(t)=\alpha(t)J(t)$, we have $\dot{\varphi}(t)x(t)+\beta(t) \cos \varphi(t)=-2\alpha(t)\beta(t)x(t)$.
It follows that 
\begin{align*}
\dot{\varphi}(t) \cos \varphi(t) x(t)+\beta(t)\cos \varphi^2(t)&=-2\alpha(t)\beta(t) x(t) \cos \varphi(t),\\
\dot{\varphi}(t) \sin \varphi(t) x(t)+\beta(t)\cos \varphi(t) \sin \varphi(t)&=-2\alpha(t)\beta(t)x(t) \sin \varphi(t).
\end{align*}
We define $X(t)=x(t)\sin \varphi(t)+i x(t)\cos \varphi(t)$, where $i$ is the imaginary unit. 
Then 
$$
\dot{X}(t)=\dot{\varphi}(t)x(t) \cos \varphi(t)-\beta(t)(1-\cos^2\varphi(t))-i (\dot{\varphi}(t)x(t)+\beta(t)\sin \varphi(t)).
$$
By a direct calculation, we have $\dot{X}(t)-2i\alpha(t)\beta(t)X(t)=-\beta(t)$.
A solution of the first order ordinary linear differential equation is given by $X(t)=(F(t)-i G(t))(\cos \eta(t)+i\sin \eta(t))$, where 
$$
F(t)=-\int \beta(t)\cos \eta(t) dt, \ G(t)=-\int \beta(t)\sin \eta(t) dt, \  \eta(t)=2 \int \alpha(t)\beta(t) dt.
$$
It follows that 
\begin{align*}
x(t) \cos \varphi(t)&=F(t)\sin \eta(t)-G(t)\cos \eta(t), \\
x(t) \sin \varphi(t)&=F(t)\cos \eta(t) +G(t) \sin \eta(t).
\end{align*}
Then we have $x^2(t)=F(t)^2+G(t)^2$.
Since $x(t_0)>0$, we have $x(t)=(F(t)^2+G(t)^2)^\frac{1}{2}$, $\sin \varphi(t)$ and $\cos \varphi(t)$ around $t_0$.
Moreover, by $\dot{z}(t)=\beta(t) \cos \varphi(t)$, we have 
$$
z(t) = \int \frac{\beta(t)}{x(t)} \left( F(t)\sin \eta(t)-G(t)\cos \eta(t)\right) dt.
$$
\enD
\begin{remark}\label{rmk:mean}{\rm
Let $(\gamma,\nu)$ be a Legendre curve constructed by Theorem \ref{mean-curvature}. 
Then the curvature $(\dot{\varphi},\beta)$ of $(\gamma,\nu)$ is calculated as 
\begin{equation}\label{eq:mean-dotphi}
(\dot{\varphi}(t),\beta(t))=(-\beta(t)X(t),\beta(t))\quad \left(X(t)=\dfrac{F(t)\sin\eta(t)-G(t)\cos\eta(t)}{F(t)^2+G(t)^2}+2\alpha(t)\right).
\end{equation} 
Thus the curvature $(\dot{\varphi}(t),\beta(t))$ as in \eqref{eq:mean-dotphi} vanishes at a singular point $t_0$ of $\gamma$. 
This implies that the curve $\gamma$ constructed by \eqref{mean-xz} is a frontal, but not a front at $t_0$ because $(\dot{\varphi}(t_0),\beta(t_0))\neq(0,0)$ 
when $(\gamma,\nu)$ is a Legendre immersion (see \cite{Fukunaga-Takahashi1}). 
Therefore the surface of revolution $\bz$ of $\gamma$ is a frontal but not a front by Proposition \ref{front-frontal}.
For the case of $\gamma$ being a front, see \cite{Matins-Saji-Santos-Teramoto}.}
\end{remark}

\begin{theorem}\label{mean-singular}
{Let $(\gamma,\nu):(I,t_0) \to \R^2 \times S^1$ be a Legendre curve of the form $\gamma(t)=(x(t),z(t))$ and $\nu(t)=(\cos \varphi(t), \sin \varphi(t))$} 
with the curvature $(\dot{\varphi},\beta)$. 
\par
$(1)$ Suppose that $x(t_0)>0$ and given a smooth function $J:(I,t_0) \to \R$. 
Then  $\gamma$ is given by 
\begin{equation*}
x(t) = \left(2 \int J(t) \sin \varphi(t) dt \right)^\frac{1}{2}, \ 
z(t) = - \int \frac{J(t)}{x(t)} \cos \varphi(t) dt.
\end{equation*}
\par
$(2)$ Suppose that $\cos \varphi(t_0) \not=0$ and given a smooth function  $H: (I,t_0) \to \R$. 
Then  $\gamma$ is given by 
\begin{equation*}
x(t) = \frac{1}{\cos \varphi(t)} \left(-\int 2H(t) \sin\varphi(t) dt \right), \ 
z(t) = \int \left(2H(t)-\dot{\varphi}(t) x(t) \right) dt.
\end{equation*}
\end{theorem}
\demo
$(1)$ By the assumption, we have $J(t)=-\beta(t)x(t)$. 
Since {the differential $\dot{x}(t)$ of $x(t)$ by $t$ is} $\dot{x}(t)=-\beta(t)\sin \varphi(t)$, we have
$$
\frac{d}{dt} x^2(t)=2x(t)\dot{x}(t)=-2x(t)\beta(t)\sin \varphi(t)=2J(t)\sin \varphi(t).
$$
Hence, 
$$
x^2(t)=2 \int J(t) \sin \varphi(t) dt.
$$ 
By $x(t_0)>0$, we have 
$$
x(t)=\left(2\int^t_{t_0} J(t) \sin \varphi(t) dt\right)^{\frac{1}{2}}.
$$
Since $\dot{z}(t)=\beta(t)\cos \varphi(t)$ and $\beta(t)=-J(t)/x(t)$, we have 
$$
z(t) = -\int \frac{J(t)}{x(t)} \cos \varphi(t) dt.
$$
\par
$(2)$ By the assumption, we have $H(t)=(1/2)(\dot{\varphi}(t)x(t)+\beta(t)\cos \varphi(t))$. 
Since $\dot{x}(t)=-\beta(t)\sin \varphi(t)$, we have
\begin{equation*}
(\cos \varphi(t))\dot{x}(t)=-\beta(t)\cos \varphi(t)\sin \varphi(t)
=(\dot{\varphi}(t)\sin \varphi(t)) x(t)-2H(t)\sin \varphi(t).
\end{equation*}
By $\cos \varphi(t_0) \not=0$, a solution of the first order ordinary linear differential equation 
$(\cos \varphi(t))\dot{x}(t)-(\dot{\varphi}(t)\sin \varphi(t)) x(t)=-2H(t)\sin \varphi(t)$ is given by 
$$
x(t) = \frac{1}{\cos \varphi(t)} \left(-\int 2H(t) \sin\varphi(t) dt \right).
$$
Since $\dot{z}(t)=\beta(t)\cos \varphi(t)$ and $\beta(t)=(2H(t)-\dot{\varphi}(t)x(t))/\cos \varphi(t)$, we have 
$$
z(t) = \int \left(2H(t)-\dot{\varphi}(t) x(t) \right) dt.
$$
\enD

We have already known that the surface of revolution of a parallel curve of a Legendre curve is a parallel surface of the surface of revolution of the frontal (cf. \cite{Gray}).

\begin{proposition}\label{parallel}
Let $(\gamma,\nu):I \to \R^2 \times S^1$ be a Legendre curve.
The parallel surface $\bz^\lambda:I \times [0,2\pi) \to \R^3$ of the surface of revolution of the frontal $\gamma:I \to \R^2$ is the surface of revolution of a parallel curve $\gamma^\lambda:I \to \R^2$ of the Legendre curve $(\gamma,\nu)$.
\end{proposition}
\demo
By definition, the parallel surface of $(\bz,\bn^z,\bs^z):I \times [0,2\pi) \to \R^3 \times \Delta$ is given by $\bz^\lambda:I \times [0,2\pi) \to \R^3,$
$$
\bz^\lambda(t,\theta)=\bz(t,\theta)+\lambda \bn^z(t,\theta)=
((x(t)+\lambda a(t)) \cos \theta,(x(t)+\lambda a(t))\sin \theta, z(t)+\lambda a(t)).
$$ 
\par
On the other hand, the parallel curve of the Legendre curve $(\gamma,\nu):I \to \R^2 \times S^1$ is given by $\gamma^\lambda:I \to \R^2, \gamma^\lambda(t)=\gamma(t)+\lambda \nu(t)=(x(t)+\lambda a(t),z(t)+\lambda b(t))$. 
It follows that the surface of revolution of the parallel curve $\gamma^\lambda$ is given by $\bz^\lambda:I \times [0,2\pi) \to \R^3$.
\enD

Let $(\gamma,\nu):I \to \R^2 \times S^1$ be a Legendre curve and $(\bz,\bn^z,\bs^z):I \times [0,2\pi) \to \R^3 \times \Delta$ be the surface of revolution. 
We consider the evolute $\mathcal{E}v(\bz)$ of the surface of revolution $\bz$ of the frontal $\gamma$.
In the following, we denote by $\mathcal{E}v^z$ the evolute of $\bz$ instead of $\mathcal{E}v(\bz)$.

\begin{proposition}\label{evolute}
Let $(\gamma,\nu):I \to \R^2 \times S^1$ be a Legendre curve with the curvature $(\ell,\beta)$.
\par
$(1)$ Suppose that $K^z_F(t,\theta) \not=0$ for all $(t,\theta) \in I \times [0,2\pi)$.
One of the evolute of the surface of revolution of the front is the surface of revolution of the evolute of the front.
The other evolute of the surface of revolution of the front is given by 
$\mathcal{E}v^z(t,\theta)=(0,0,z(t)-x(t)b(t)/a(t))$.
\par
$(2)$ Suppose that $\ell(t) \not=0$ for all $t \in I$.
At least an evolute of the surface of revolution of the front is the surface of revolution of the evolute of the front.
\par
$(3)$ Suppose that $a(t) \not=0$ for all $t \in I$.
At least an evolute of the surface of revolution of the front is given by $\mathcal{E}v^z(t,\theta)=(0,0,z(t)-x(t)b(t)/a(t))$.
\end{proposition}
\demo
$(1)$ Since $K^z_F(t,\theta)=-a(t)\ell(t) \not=0$, we have $a(t) \not=0$ and $\ell(t) \not=0$ for all $t \in I$. 
Since 
{
\begin{align*}
K^z_F(t,\theta)\lambda^2-2H^z_F(t,\theta)\lambda+J^z_F(t,\theta)
&=-a(t)\ell(t)\lambda^2-(x(t)\ell(t)+\beta(t)a(t))\lambda-x(t)\beta(t)\\
&=-(a(t)\lambda+x(t))(\ell(t)\lambda+\beta(t)),
\end{align*}
}
the solutions of the equation $K^z_F(t,\theta)\lambda^2-2H^z_F(t,\theta)\lambda+J^z_F(t,\theta)=0$ are given by $\lambda=-\beta(t)/\ell(t)$ and $\lambda=-x(t)/a(t)$.
Therefore, one of the evolute is given by 
\begin{align*}
\mathcal{E}v^z(t,\theta)&=\bz(t,\theta)-\frac{\beta(t)}{\ell(t)}\bn^z(t,\theta)\\
&=\left((x(t) \cos \theta,x(t)\sin \theta,z(t)\right)-\frac{\beta(t)}{\ell(t)} \left(a(t) \cos \theta,a(t)\sin \theta,b(t)\right)\\
&=\left( \left(x(t)-\frac{\beta(t)}{\ell(t)}a(t)\right)\cos \theta, \left(x(t)-\frac{\beta(t)}{\ell(t)}a(t)\right)\sin \theta, z(t)-\frac{\beta(t)}{\ell(t)}b(t)\right).
\end{align*}
Hence it is the surface of revolution of the evolute $\mathcal{E}v(\gamma)(t)=\gamma(t)-(\beta(t)/\ell(t))\nu(t)$.
\par
The other evolute is given by 
\begin{align*}
\mathcal{E}v^z(t,\theta)&=\bz(t,\theta)-\frac{x(t)}{a(t)}\bn^z(t,\theta)\\
&=\left((x(t) \cos \theta,x(t)\sin \theta,z(t)\right)-\frac{x(t)}{a(t)} \left(a(t) \cos \theta,a(t)\sin \theta,b(t)\right)\\
&=\left(0,0,z(t)-\frac{x(t)}{a(t)}b(t)\right).
\end{align*}
\par
$(2)$ Suppose that $\ell(t) \not=0$ for all $t \in I$. 
Then we have at least a solution $\lambda=-\beta(t)/\ell(t)$ of the equation $K^z_F(t,\theta)\lambda^2-2H^z_F(t,\theta)\lambda+J^z_F(t,\theta)=0$. 
Therefore, by the same calculation of $(1)$, 
at least an evolute of the surface of revolution of the front is the surface of revolution of evolute of the front. 
\par
$(3)$ By the same argument as in $(2)$, we have the result.
\enD
\begin{remark}\label{rem:evolute2}{\rm
If $\ell(t)=0$, $a(t) \not=0$ and $\beta(t) \not=0$ for all $t \in I$, then  we have the evolute, see Remark \ref{rem:evolute1} and Proposition \ref{evolute} $(3)$. 
Moreover, even if $\beta(t)=0$ for isolated points, we also define an evolute of the surface of revolution of the frontal as the same by the continuous property. 
}
\end{remark}

We consider classification problems. 
{We use notations $x_0=x(t_0)$, $z_0=z(t_0)$ and $\bz_0=\bz(t_0,\theta_0)=(x_0\cos \theta_0,x_0\sin \theta_0,z_0)$ for 
smooth curve $\gamma=(x,z):(I,t_0)\to(\R^2,(x_0,z_0))$ and its surface of revolution $\bz$ around $z$-axis, respectively.
We also use a notation $\gamma_0=\gamma(t_0)=(x_0,z_0)$.}

First we consider the case of $x_0>0$.

\begin{theorem}\label{diffeo-positive}
Let $\gamma:(I,t_0) \to (\R^2,(x_0,z_0))$ and $\widetilde{\gamma}:(\widetilde{I},\widetilde{t}_0) \to (\R^2,(\widetilde{x}_0,\widetilde{z}_0))$ be smooth curves 
with $x_0, \widetilde{x}_0>0$, 
let {$\bz:(I \times [0,2\pi),(t_0,\theta_0)) \to (\R^2 \setminus \{0\} \times \R,\bz_0)$ 
and $\widetilde{\bz}:(\widetilde{I} \times [0,2\pi),(\widetilde{t}_0,\theta_0)) \to (\R^2 \setminus \{0\} \times \R,\widetilde{\bz}_0)$} 
be surfaces of revolution of the frontals around $z$-axis, respectively. 
\par
$(1)$ If there exist diffeomorphism germs $\phi:(I,t_0) \to (\widetilde{I},\widetilde{t_0})$ 
and {$\psi: (\R^2,\gamma_0)\to (\R^2,\widetilde{\gamma}_0)$, $\psi(X,Z)=(\psi_1(X,Z),\psi_2(X,Z))$} such that $\psi \circ \gamma=\widetilde{\gamma} \circ \phi$, 
then there exist diffeomorphism germs $\Phi:(I \times [0,2\pi),(t_0,\theta_0)) \to (\widetilde{I} \times [0,2\pi),(\widetilde{t_0},\theta_0))$ of the form $\Phi(t,\theta)=(\phi(t),\theta)$ 
and {$\Psi:(\R^2 \setminus \{0\} \times \R,\bz_0) \to  (\R^2 \setminus \{0\} \times \R,\widetilde{\bz}_0)$} of the form 
$$
\Psi(X,Y,Z)=\left(\frac{X\psi_1(\sqrt{X^2+Y^2},Z)}{\sqrt{X^2+Y^2}}, \frac{Y\psi_1(\sqrt{X^2+Y^2},Z)}{\sqrt{X^2+Y^2}}, \psi_2(\sqrt{X^2+Y^2},Z)\right)
$$ 
such that $\Psi \circ \bz=\widetilde{\bz} \circ \Phi$.
\par
$(2)$ If there exist diffeomorphism germs $\Phi:(I \times [0,2\pi),(t_0,\theta_0)) \to (\widetilde{I} \times [0,2\pi),(\widetilde{t_0},\theta_0))$ of the form $\Phi(t,\theta)=(\phi(t),\theta)$ 
and {$\Psi:(\R^2 \setminus \{0\} \times \R,\bz_0) \to  (\R^2 \setminus \{0\} \times \R,\widetilde{\bz}_0)$, $\Psi(X,Y,Z)=(\Psi_1(X,Y,Z),\Psi_2(X,Y,Z),\Psi_3(X,Y,Z))$} 
such that $\Psi \circ \bz=\widetilde{\bz} \circ \Phi$, 
then there exists a diffeomorphism germ {$\psi: (\R^2,\gamma_0)\to (\R^2,\widetilde{\gamma}_0)$}, $\psi(X,Z)=(\psi_1(X,Z),\psi_2(X,Z))
$ of the form 
\begin{align*}
\psi_1(X,Z)&= \Psi_1(X\cos \theta_0,X\sin \theta_0,Z)\cos \theta_0+\Psi_2(X\cos \theta_0,X\sin \theta_0,Z)\sin \theta_0,\\
\psi_2(X,Z)&=\Psi_3(X\cos \theta_0,X\sin \theta_0,Z)
\end{align*}
such that $\psi \circ \gamma=\widetilde{\gamma} \circ \phi$
\end{theorem}
\demo
$(1)$ Since $\phi$ is a diffeomorphism germ, $\Phi$ is a diffeomorphism germ. 
We show that $\Psi(X,Y,Z)=(\Psi_1(X,Y,Z),\Psi_2(X,Y,Z),\Psi_3(X,Y,Z))$ is a diffeomorphism germ. 
By a direct calculation, we have 
{
\begin{align*}
\Psi_{1X}(\bz_0)&=\frac{\psi_1(x_0,z_0)\sin^2 \theta_0+x_0 \psi_{1X}(x_0,z_0) \cos^2 \theta_0}{x_0},\\
\Psi_{1Y}(\bz_0)&=\frac{(x_0\psi_{1X}(x_0,z_0)- \psi_1(x_0,z_0))\cos \theta_0 \sin \theta_0}{x_0},\\
\Psi_{1Z}(\bz_0)&=\psi_{1Z}(x_0,z_0)\cos \theta_0,\\
\Psi_{2X}(\bz_0)&=\frac{(x_0\psi_{1X}(x_0,z_0)- \psi_1(x_0,z_0))\cos \theta_0 \sin \theta_0}{x_0},\\
\Psi_{2Y}(\bz_0)&=\frac{\psi_1(x_0,z_0)\cos^2 \theta_0+x_0 \psi_{1X}(x_0,z_0) \sin^2 \theta_0}{x_0},\\
\Psi_{2Z}(\bz_0)&=\psi_{1Z}(x_0,z_0)\sin \theta_0, \\
\Psi_{3X}(\bz_0)&=\psi_{2X}(x_0,z_0)\cos \theta_0, \\
\Psi_{3Y}(\bz_0)&=\psi_{2X}(x_0,z_0)\sin \theta_0, \\
\Psi_{3Z}(\bz_0)&=\psi_{2Z}(x_0,z_0).
\end{align*}
}
Then we can show that the determinant of the Jacobi matrix of $\Psi$ at {$\bz_0$ is non-zero}. 
Hence $\Psi$ is a diffeomorphism germ.
By the assumption, we have 
$$
(\widetilde{x}(\phi(t)),\widetilde{z}(\phi(t)))=(\psi_1(x(t),z(t)),\psi_2(x(t),z(t))).
$$
It follows that 
\begin{align*}
\Psi \circ \bz(t,\theta)&=\Psi (x(t)\cos \theta,x(t)\sin \theta,z(t))\\
&=(\psi_1(x(t),z(t)) \cos \theta ,\psi_1(x(t),z(t)) \sin \theta,\psi_2(x(t),z(t)))\\
&=(\widetilde{x}(\phi(t))\cos \theta, \widetilde{x}(\phi(t))\sin \theta,\widetilde{z}(\phi(t)))\\
&=\widetilde{\bz} \circ \Phi(t,\theta).
\end{align*}
$(2)$ 
Since $\Phi$ is a diffeomorphism germ, $\phi$ is a diffeomorphism germ. 
We show that $\psi(X,Z)=(\psi_1(X,Z),\psi_2(X,Z))$ is a diffeomorphism germ. 
By a direct calculation, we have 
{
\begin{align*}
\psi_{1X}(x_0,z_0)&=\Psi_{1X}(\bz_0)\cos^2 \theta_0+\Psi_{1Y}(\bz_0)\cos \theta_0 \sin \theta_0\\
&\quad+\Psi_{2X}(\bz_0)\cos \theta_0 \sin \theta_0+\Psi_{2Y}(\bz_0) \sin^2 \theta_0,\\
\psi_{1Z}(x_0,z_0)&=\Psi_{1Z}(\bz_0)\cos \theta_0+\Psi_{2Z}(\bz_0)\sin \theta_0,\\
\psi_{2X}(x_0,z_0)&=\Psi_{3X}(\bz_0)\cos \theta_0+\Psi_{3Y}(\bz_0)\sin \theta_0,\\
\psi_{2Z}(x_0,z_0)&=\Psi_{3Z}(\bz_0).
\end{align*}
}
Since $\Psi \circ \bz=\widetilde{\bz} \circ \Phi$, we have 
\begin{align*}
\widetilde{x}(\phi(t)) \cos \theta &= {\Psi}_1(x(t)\cos \theta,x(t) \sin \theta,z(t)),\\
\widetilde{x}(\phi(t)) \sin \theta &= {\Psi}_2(x(t)\cos \theta,x(t) \sin \theta,z(t)),\\
\widetilde{z}(\phi(t)) &= {\Psi}_3(x(t)\cos \theta,x(t) \sin \theta,z(t))
\end{align*}
for all $(t, \theta)$. 
By differentiating with respect to $\theta$, we have 
\begin{align*}
-\widetilde{x}(\phi(t)) \sin \theta &= {\Psi}_{1X}(x(t)\cos \theta,x(t) \sin \theta,z(t))(-x(t)\sin \theta)\\
&\quad+{\Psi}_{1Y}(x(t)\cos \theta,x(t) \sin \theta,z(t))x(t)\cos \theta \\
&=-{\Psi}_2(x(t)\cos \theta,x(t) \sin \theta,z(t)),\\
\widetilde{x}(\phi(t)) \cos \theta &= {\Psi}_{2X}(x(t)\cos \theta,x(t) \sin \theta,z(t))(-x(t)\sin \theta)\\
&\quad+{\Psi}_{2Y}(x(t)\cos \theta,x(t) \sin \theta,z(t))x(t)\cos \theta\\
&={\Psi}_1(x(t)\cos \theta,x(t) \sin \theta,z(t)),\\
0 &= {\Psi}_{3X}(x(t)\cos \theta,x(t) \sin \theta,z(t))(-x(t) \sin \theta)\\
&\quad+{\Psi}_{3Y}(x(t)\cos \theta,x(t) \sin \theta,z(t))x(t) \cos \theta. 
\end{align*}
By using the above, we can show that the determinant of the Jacobi matrix of $\psi$ at {$\gamma_0=(x_0,z_0)$} is non-zero.
Hence $\psi$ is a diffeomorphism germ.
Moreover, 
{
\begin{align*}
\psi \circ \gamma(t)&=(\psi_1(x(t),z(t)),\psi_2(x(t),z(t)))\\
&=(\Psi_1(x(t)\cos \theta_0,x(t)\sin \theta_0,z(t))\cos \theta_0\\
&+\Psi_2(x(t)\cos \theta_0,x(t)\sin \theta_0,z(t))\sin \theta_0, 
\Psi_3(x(t)\cos \theta_0,x(t)\sin \theta_0,z(t)))\\
&=(\widetilde{x}(\phi(t))\cos^2 \theta_0+\widetilde{x}(\phi(t))\sin^2 \theta_0, \widetilde{z}(\phi(t)))\\
&=\widetilde{\gamma} \circ \phi(t).
\end{align*}
}
This completes the proof of {Theorem \ref{diffeo-positive}}.
\enD

To characterize singularities of surfaces of revolution applying Theorem \ref{diffeo-positive}, 
we give some definitions (cf. \cite{Izumiya-book}).
\begin{definition}{\rm
$(1)$ Let $f$ and $g:(\R^m,0)\to(\R^n,0)$ be smooth map-germs. 
Then $f$ is {\it $\mathcal{A}$-equivalent} to $g$ if there exist diffeomorphism germs $\phi:(\R^m,0)\to(\R^m,0)$ and $\Phi:(\R^n,0)\to(\R^n,0)$ 
such that $g=\Phi\circ f\circ \phi^{-1}$ holds. 
If the diffeomorphism germ $\Phi$ (respectively, $\phi$) appeared in above is the identity map, 
we say that $f$ is {\it $\mathcal{R}$-equivalent} (respectively, {\it $\mathcal{L}$-equivalent}) to $g$.
\par
$(2)$ Let $\gamma:(I,t_0)\to(\R^2,0)$ be a smooth curve. 
We say that $\gamma$ at $t_0$ is a {\it $j/i$-cusp}, where $(i,j)=(2,3),(2,5),(3,4),(3,5)$ if $\gamma$ is $\mathcal{A}$-equivalent to the germ $t\mapsto(t^i,t^j)$ at the origin.
\par
$(3)$ Let $f:(\R^2,0)\to(\R^3,0)$ be a smooth map. 
We say that $f$ at $0$ is a {\it $j/i$-cuspidal edge}, where $(i,j)=(2,3),(2,5),(3,4),(3,5)$ if $f$ is $\mathcal{A}$-equivalent to the germ $(u,v)\mapsto(u,v^i,v^j)$ at the origin.}
\end{definition}

We note that curves with $j/i$-cusps are frontal (curves). 
Moreover, surfaces with $j/i$-cuspidal edges are not only frontal (surfaces), 
but also framed base surfaces (see \cite{Fukunaga-Takahashi1,Fukunaga-Takahashi2,Fukunaga-Takahashi3}).
{In particular, $3/2$-cusps and $4/3$-cusps (respectively, $3/2$-cuspidal edges and $4/3$-cuspidal edges) are front singularities.} 
As a corollary of Theorem \ref{diffeo-positive}, we have the following.

\begin{corollary}
Let $\gamma=(x,z):(I,t_0)\to(\R^2,(x_0,z_0))$ be a smooth curve with $x_0>0$. 
Then $\gamma$ at $t_0$ is a $j/i$-cusp 
if and only if the surface of revolution $\bz(t,\theta)$ of $\gamma$ 
around the $z$-axis at $(t_0,\theta_0)$ is a $j/i$-cuspidal edge for any $\theta_0$.
\end{corollary}

For curves with $j/i$-cusps, the following criteria are known (cf. \cite{Bruce-Gaffney, Porteous}). 
\begin{proposition}\label{fact:criteria}
Let $\gamma:I\to\R^2$ be a smooth curve and $t_0\in I$ a singular point of $\gamma$, 
namely, $\dot{\gamma}(t_0)=0$. 
Then the following assertions hold.
\par
$(1)$ $\gamma$ has a $3/2$-cusp at $t_0$ if and only if $\det(\ddot{\gamma},\dddot{\gamma})(t_0)\neq0$.
\par
$(2)$ $\gamma$ has a $5/2$-cusp at $t_0$ if and only if 
$\ddot{\gamma}(t_0)\neq0$, $\dddot{\gamma}(t_0)=C\ddot{\gamma}(t_0)$ for some constant $C\in\R$ 
and $\det(\ddot{\gamma},3\gamma^{(5)}-10C\gamma^{(4)})(t_0)\neq0$.
\par
$(3)$ $\gamma$ has a $4/3$-cusp at $t_0$ if and only if $\ddot{\gamma}(t_0)=0$ and 
$\det(\dddot{\gamma},\gamma^{(4)})(t_0)\neq0$.
\par
$(4)$ $\gamma$ has a $5/3$-cusp at $t_0$ if and only if $\ddot{\gamma}(t_0)=0$, 
$\det(\dddot{\gamma},\gamma^{(4)})(t_0)=0$ and $\det(\dddot{\gamma},\gamma^{(5)})(t_0)\neq0$.
\end{proposition}
Using Proposition \ref{fact:criteria}, we show the following.
\begin{theorem}\label{thm:criteria}
{
Let $(\gamma,\nu):(I,t_0)\to\R^2\times S^1$ 
be a Legendre curve of the form $\gamma(t)=(x(t),z(t))$, $\nu(t)=(\cos\varphi(t),\sin\varphi(t))$ with the curvature $(\dot{\varphi},\beta)$.} 
Assume that $\beta(t_0)=0$. 
Then we have the following.
\par
$(1)$ $\gamma$ is a $3/2$-cusp at $t_0$ if and only if $\dot{\beta}(t_0)\dot{\varphi}(t_0)\neq0$.
\par
$(2)$ $\gamma$ is a $5/2$-cusp at $t_0$ if and only if $\dot{\beta}(t_0)\neq0$, $\dot{\varphi}(t_0)=0$ and 
$\ddot{\beta}(t_0)\ddot{\varphi}(t_0)-\dot{\beta}(t_0)\dddot{\varphi}(t_0)\neq0$.
\par
$(3)$ $\gamma$ is a $4/3$-cusp at $t_0$ if and only if $\dot{\beta}(t_0)=0$ and $\ddot{\beta}(t_0)\dot{\varphi}(t_0)\neq0$.
\par
$(4)$ $\gamma$ is a $5/3$-cusp at $t_0$ if and only if $\dot{\beta}(t_0)=\dot{\varphi}(t_0)=0$ and $\ddot{\beta}(t_0)\ddot{\varphi}(t_0)\neq0$.
\end{theorem}
\demo
By the Frenet type formula \eqref{Frenet.frontal}, we have $\dot{\gamma}(t)=\beta(t)\bmu(t)$, $\dot{\nu}(t)=\dot{\varphi}(t)\bmu(t)$ and $\dot{\bmu}(t)=-\dot{\varphi}(t)\nu(t)$, where $\bmu(t)=(-\sin\varphi(t),\cos\varphi(t))$. 
We consider differentials of $\gamma(t)$. 
By direct calculations, we see that 
\begin{align*}
\ddot{\gamma}(t)&=\dot{\beta}(t)\bmu(t)-\beta(t)\dot{\varphi}(t)\nu(t),\\
\dddot{\gamma}(t)&=(\ddot{\beta}(t)-\beta(t)(\dot{\varphi})^2(t)\bmu(t)-(2\dot{\beta}(t)\dot{\varphi}(t)+\beta(t)\ddot{\varphi}(t))\nu(t),\\
\gamma^{(4)}(t)&=(\dddot{\beta}(t)-3\dot{\beta}(t)(\dot{\varphi})^2(t)-3\beta(t)\dot{\varphi}(t)\ddot{\varphi}(t))\bmu(t)\\
&\quad -(3\ddot{\beta}(t)\dot{\varphi}(t)+3\dot{\beta}(t)\ddot{\varphi}(t)+\beta(t)\dddot{\varphi}(t)-\beta(t)(\dot{\varphi})^3(t)\nu(t).
\end{align*}
Since $\beta(t_0)=0$, we have $\det(\ddot{\gamma},\dddot{\gamma})(t_0)=2(\dot{\beta})^2(t_0)\dot{\varphi}(t_0)$. 
Thus the assertion $(1)$ holds by Proposition \ref{fact:criteria} $(1)$.
 
We show $(3)$ by using Proposition \ref{fact:criteria} $(3)$. 
The condition for $\ddot{\gamma}(t_0)=0$ is $\dot{\beta}(t_0)=0$. 
In this case, $\dddot{\gamma}(t_0)=\ddot{\beta}(t_0)\bmu(t_0)$, 
and hence $\dddot{\gamma}(t_0)\neq0$ if and only if $\ddot{\beta}(t_0)\neq0$. 
Moreover, $\det(\dddot{\gamma},\gamma^{(4)})(t_0)=3(\ddot{\beta})^2(t_0)\dot{\varphi}(t_0)$ holds. 
Thus we have $(3)$.  

We next consider $(2)$. 
Since $\ddot{\gamma}(t_0)\neq0$, we note that $\dot{\beta}(t_0)\neq0$ holds. 
Two vectors $\ddot{\gamma}(t_0)$ and $\dddot{\gamma}(t_0)$ are parallel 
if and only if $\det(\ddot{\gamma},\dddot{\gamma})(t_0)=0$. 
This is equivalent to $\dot{\varphi}(t_0)=0$, 
and hence we see that $\dddot{\gamma}(t_0)=C\ddot{\gamma}(t_0)$, where $C=\ddot{\beta}(t_0)/\dot{\beta}(t_0)$. 
Further, $\gamma^{(4)}(t_0)$ and $\gamma^{(5)}(t_0)$ are calculated as
\begin{align*}
\gamma^{(4)}(t_0)&=\dddot{\beta}(t_0)\bmu(t_0)-3\dot{\beta}(t_0)\ddot{\varphi}(t_0)\nu(t_0),\\
\gamma^{(5)}(t_0)&=\beta^{(4)}(t_0)\bmu(t_0)-(6\ddot{\beta}(t_0)\ddot{\varphi}(t_0)+4\dot{\beta}(t_0)\dddot{\varphi}(t_0))\nu(t_0).
\end{align*}
Thus 
$$\det\left(\ddot{\gamma},3\gamma^{(5)}-10C\gamma^{(4)}\right)(t_0)
=-12\dot{\beta}(t_0)(\ddot{\beta}(t_0)\ddot{\varphi}(t_0)-\dot{\beta}(t_0)\dddot{\varphi}(t_0))$$
holds. 
By Proposition \ref{fact:criteria} $(2)$, we see that the assertion $(2)$ holds.

Finally we show $(4)$. 
Since $\ddot{\gamma}(t_0)=\det(\dddot{\gamma},\gamma^{(4)})(t_0)=0$ and $\dddot{\gamma}(t_0)\neq0$, 
we see that $\dot{\beta}(t_0)=\dot{\varphi}(t_0)=0$ and $\ddot{\beta}(t_0)\neq0$. 
Under these conditions, $\gamma^{(5)}(t_0)$ is given by 
{
$$\gamma^{(5)}(t_0)=\beta^{(4)}(t_0)\bmu(t_0)-6\ddot{\beta}(t_0)\ddot{\varphi}(t_0)\nu(t_0).$$
}
Hence we have 
{
$$\det(\dddot{\gamma},\gamma^{(5)})(t_0)=6(\ddot{\beta})^2(t_0)\ddot{\varphi}(t_0).$$
}
Thus we have the assertion by Proposition \ref{fact:criteria} $(4)$.
\enD
If we use the notation for a Legendre curve $(\gamma,\nu)$ with the curvature $(\ell,\beta)$, 
assertions in Theorem \ref{thm:criteria} also hold by replacing $\dot{\varphi}$ with $\ell$. 

Let $f:(I,t_0)\to\R$ be a function germ and $k$ be a non-negative integer. 
Then $f$ has a {\it zero of order $(k+1)$} at $t_0$ if 
\begin{equation}\label{eq:order}
f(t_0)=\dot{f}(t_0)=\ddot{f}(t_0)=\cdots=f^{(k)}(t_0)=0,\quad f^{(k+1)}(t_0)\neq0,
\end{equation}
where $\dot{f}=df/dt$ and $f^{(n)}=d^nf/dt^n$ for a positive integer $n$. 
In this case, we write $\ord(f)(t_0)=k+1$. 
For a Legendre curve $(\gamma,\nu)$ with the curvature $(\dot{\varphi},\beta)$, 
we remark that $\beta$ satisfies $\ord(\beta)(t_0)\geq1$ at a singular point $t_0$ of $\gamma$ since $\beta(t_0)=0$.

We now use notations $\nu(t)=(a(t),b(t))$ for the unit normal and $(\ell(t),\beta(t))$ for the curvature 
instead of $\nu(t)=(\cos\varphi(t),\sin\varphi(t))$ and $(\dot{\varphi}(t),\beta(t))$, respectively. 
\begin{proposition}\label{prop:Gauss-front}
Let $(\gamma,\nu):I\to\R^2\times S^1$, $\gamma(t)=(x(t),z(t))$, $\nu(t)=(a(t),b(t))$ be a Legendre curve which is given by 
$(1)$ in Theorem \ref{Gauss-curvature}. 
Suppose that $x(t_0)\neq0$ and $\alpha(t_0)\neq0$ at a singular point $t_0$ of $\gamma$. 
Then we have the following.
\par
$(1)$ If $a(t_0)\neq0$, then $\gamma$ is a frontal but not a front at $t_0$.
\par
$(2)$ Suppose that $\ord(a)(t_0)=m+1$ and $\ord(\beta)(t_0)=n+1$. 
Then $\gamma$ is a front at $t_0$ if and only if $\ord(a)(t_0)=\ord(\beta)(t_0)$.
\end{proposition}
\demo
Without loss of generality, we may assume that $t_0=0$. 
By the proof of $(1)$ in Theorem \ref{Gauss-curvature}, the function $\ell(t)$ satisfies
$$\ell(t)a(t)=\alpha(t)\beta(t)x(t).$$
If $a(0)\neq0$, the function $\ell$ can be given as 
$$\ell(t)=\dfrac{\alpha(t)\beta(t)x(t)}{a(t)}$$
around $0$. 
Since $\beta(0)=0$, we have $(\ell(0),\beta(0))=(0,0)$. 
Hence, $\gamma$ is a frontal but not a front at $t=0$. 
Thus the assertion $(1)$ holds.

We show the assertion $(2)$. 
We now assume that $\ord(a)(0)=m+1$ and $\ord(\beta)(0)=n+1$.  
Then by the division lemma, there exist smooth functions $\hat{a}(t)$ and $\hat{\beta}(t)$ around $t=0$ such that 
$a(t)=t^{m+1}\hat{a}(t)$ and $\beta(t)=t^{n+1}\hat{\beta}(t)$.
We note that $\hat{a}(0)\neq0$ and $\hat{\beta}(0)\neq0$. 
If $m>n$, then we have $t^{m-n} \hat{a}(t) \ell(t)=\alpha(t) \hat{\beta}(t) x(t)$. 
It follows that $\alpha(0) \hat{\beta}(0) x(0)=0$. 
This contradicts the fact that $\alpha(0) \hat{\beta}(0) x(0) \not=0$.
Therefore, $m\leq n$. 
Then the function $\ell(t)$ can be expressed as 
\begin{equation}\label{eq:Gauss-ell}
\ell(t)=\dfrac{t^{n-m}\alpha(t)\hat{\beta}(t)x(t)}{\hat{a}(t)}.
\end{equation}
Thus $\ell(0)\neq0$ if and only if $n=m$. 
This implies that $\ord(a)(0)=\ord(\beta)(0)$. 
Therefore we have the assertion $(2)$.
\enD

We have the following characterizations of singularities 
for the case of constant Gauss curvature surfaces of revolution (see $(1)$ in Theorem \ref{Gauss-curvature}). 
\begin{proposition}\label{constant-Gauss}
Under the same assumptions as in Theorem \ref{Gauss-curvature} $(1)$ with conditions that the Legendre curve $(\gamma,\nu):(I,t_0) \to \R^2 \times S^1$,  $\gamma(t)=(x(t),z(t)), \nu(t)=(a(t),b(t))$ satisfy $x(t_0) \neq 0$ and $\ord(a)(t_0)=m+1 \le \ord(\beta)(t_0)=n+1$. 
If the smooth function $\alpha:(I,t_0)\to\R$ satisfying $K(t)=\alpha(t)J(t)$ is a non-zero constant $c\in\R$, namely, $\alpha(t)=c$, 
then the curve $\gamma$ given by $(1)$ in Theorem \ref{Gauss-curvature} has
\par
$(1)$ a $3/2$-cusp at $t_0$ if and only if $\ord(a)(t_0)=\ord(\beta)(t_0)=1$,
\par
$(2)$ a $4/3$-cusp at $t_0$ if and only if $\ord(a)(t_0)=\ord(\beta)(t_0)=2$,
\par
$(3)$ a $5/3$-cusp at $t_0$ if and only if $\ord(a)(t_0)=1$ and $\ord(\beta)(t_0)=2$.

\noindent
Moreover, $\gamma$ cannot have a $5/2$-cusp at $t_0$.
\end{proposition}
\demo
By a parallel translation on the source, we may assume that $t_0=0$. 
By \eqref{eq:Gauss-ell}, the function $\ell(t)$ is given by 
$$\ell(t)=\dfrac{ct^{n-m}\hat{\beta}(t)x(t)}{\hat{a}(t)},$$
where $a(t)=t^{m+1}\hat{a}(t)$ with $\hat{a}(0)\neq0$ and $\beta(t)=t^{n+1}\hat{\beta}(t)$ with $\hat{\beta}(0)\neq0$. 
We first consider the case of $\gamma$ to be a front at $0$, that is, $m+1=\ord(a)(0)=\ord(\beta)(0)=n+1$ by Proposition \ref{prop:Gauss-front}.
Then $\dot{\beta}(0)\neq0$ if and only if $n=0$. 
By $(1)$ in Theorem \ref{thm:criteria}, we have the assertion $(1)$. 
Moreover, $\dot{\beta}(0)=0$ and $\ddot{\beta}(0)\neq0$ if and only if $n=1$, and hence we have the assertion $(2)$ 
by $(3)$ in Theorem \ref{thm:criteria}.

We next consider the case of $0\leq m<n$. 
In this case, $\gamma$ at $t=0$ is a frontal but not a front. 
Moreover, the functions $\ell$ and $\beta$ are expressed as 
$$\ell(t)=t^{n-m}A(t),\quad \beta(t)=t^{n+1}\hat{\beta}(t)\quad\left(A(t)=\dfrac{c\hat{\beta}(t)x(t)}{\hat{a}(t)}\right).$$ 
Thus $\dot{\beta}(0)\neq0$ if and only if $n=0$. 
Since $m<n$, it implies that $m<0$. 
This contradicts the fact that $m\geq0$. 
Therefore, $\gamma$ cannot have a $5/2$-cusp at $t=0$ by $(2)$ in Theorem \ref{thm:criteria}. 
We assume that $\dot{\beta}(0)=0$, that is, $n\geq1$. 
Under this situation, $\ddot{\beta}(0)\neq0$ if and only if $n=1$. 
Since $0\leq m<1=n$, we have $m=0$. 
Then it holds that $(\ell(t),\beta(t))=(tA(t),t^2\hat{\beta}(t))$. 
Noting $A(0)\neq0$ and $\hat{\beta}(0)\neq0$, this pair $(\ell(t),\beta(t))$ satisfy the conditions of $(4)$ in Theorem \ref{thm:criteria}. 
Therefore we have the assertion. 
\enD

We next consider singularities of curves obtained by Theorem \ref{mean-curvature}. 
Although for the similar case of Theorem \ref{mean-curvature}, criteria for $j/i$-cusps 
are obtained by using the data $\eta$ and $\beta$ (\cite[Proposition 2.3]{Matins-Saji-Santos-Teramoto}), 
in our setting, the following assertions hold. 
\begin{proposition}\label{constant-mean}
Under the same assumptions as in Theorem \ref{mean-curvature}, we have the following. 
\par
$(1)$ Suppose that $\beta(t_0)=0$ and the function $\alpha:(I,t_0)\to\R$ satisfying $H(t)=\alpha(t)J(t)$ is not a constant function. 
Then {$\gamma$ given by \eqref{mean-xz} does not have neither a $3/2$-cusp nor a $4/3$-cusp at $t_0$. 
Moreover, $\gamma$ has a $5/2$-cusp $($respectively, $5/3$-cusp$)$ at $t_0$
 if and only if $\dot{\beta}(t_0)\dot{\alpha}(t_0)\neq0$ $($respectively, $\dot{\beta}(t_0)=0$ and $\ddot{\beta}(t_0)\dot{\alpha}(t_0)\neq0)$}.
\par
$(2)$ Suppose that the function $\alpha:(I,t_0)\to\R$ satisfying $H(t)=\alpha(t)J(t)$ is a constant $c\in\R$, i.e., $\alpha(t)=c$. 
Then $\gamma$ given by \eqref{mean-xz} does not have $j/i$-cusps $((i,j)=(2,3),(2,5),(3,4),(3,5))$.
\end{proposition}
\demo
{$(1)$ By Remark \ref{rmk:mean}, a curve $\gamma$ obtained by Theorem \ref{mean-curvature} is a frontal but not a front. 
Thus we see that $3/2$-cusps and $4/3$-cusps do not appear on the curve $\gamma$. 
By \eqref{eq:mean-dotphi}, the curvature $(\dot{\varphi}(t),\beta(t))$ of a Legendre curve given by Theorem \ref{mean-curvature} is 
$$
\left(\dot{\varphi}(t),\beta(t)\right)
=\left(-{\beta(t)}X(t),\beta(t)\right)\quad \left(X(t)=\dfrac{F(t)\sin\eta(t)-G(t)\cos\eta(t)}{F(t)^2+G(t)^2}+2\alpha(t)\right).
$$
By the definitions of $F(t)$, $G(t)$ and $\eta(t)$, we see that $\dot{F}(t_0)=\dot{G}(t_0)=\dot{\eta}(t_0)=0$ since $\beta(t_0)=0$. 
Thus $\dot{X}(t_0)=2\dot{\alpha}(t_0)$ holds. 
On the other hand, it follows that $\ddot{\varphi}(t)=-\dot{\beta}(t)X(t)-\beta(t)\dot{X}(t)$ 
and $\dddot{\varphi}(t)=-\ddot{\beta}(t)X(t)-2\dot{\beta}(t)\dot{X}(t)-\beta(t)\ddot{X}(t)$. 
Therefore $\ddot{\varphi}(t_0)=-\dot{\beta}(t_0)X(t_0)$ and $\dddot{\varphi}(t_0)=-\ddot{\beta}(t_0)X(t_0)-2\dot{\beta}(t_0)\dot{X}(t_0)$ hold. 
By $(2)$ and $(4)$ of Theorem \ref{thm:criteria}, we have the assertion. }

$(2)$ By the proof of the assertion $(1)$ and the Theorem \ref{thm:criteria}, 
we see that the curvature $(\dot{\varphi},\beta)$ of a Legendre curve $(\gamma,\nu)$ given by Theorem \ref{mean-curvature} 
cannot satisfy the conditions for $j/i$-cusps $((i,j)=(2,3),(2,5),(3,4),(3,5))$ when $\alpha(t)=c$. 
Thus we have the conclusion.
\enD

The next, we consider the case of $x_0=0$.
\begin{theorem}\label{diffeo-zero}
Let $\gamma:(I,t_0) \to (\R^2,0)$ and $\widetilde{\gamma}:(\widetilde{I},\widetilde{t}_0) \to (\R^2,0)$ be smooth curves, 
let $\bz:(I \times [0,2\pi),(t_0,\theta_0)) \to (\R^3,0)$ and  $\widetilde{\bz}:(\widetilde{I} \times [0,2\pi),(\widetilde{t}_0,\theta_0)) \to (\R^3,0)$ 
be surfaces of revolution of the frontals  around $z$-axis, respectively. 
\par
$(1)$ If there exist diffeomorphism germs $\phi:(I,t_0) \to (\widetilde{I},\widetilde{t_0})$ and $\psi: (\R^2,0) \to (\R^2,0)$ of the form 
$\psi(X,Z)=(X,\varphi(X^2,Z))$ such that $\psi \circ \gamma=\widetilde{\gamma} \circ \phi$, then 
there exist diffeomorphism germs $\Phi:(I \times [0,2\pi),(t_0,\theta_0)) \to (\widetilde{I} \times [0,2\pi),(\widetilde{t_0},\theta_0))$ of the form $\Phi(t,\theta)=(\phi(t),\theta)$ and $\Psi:(\R^3,0) \to  (\R^3,0)$ of the form 
$\Psi(X,Y,Z)=(X,Y,\varphi (X^2+Y^2,Z))$ 
such that $\Psi \circ \bz=\widetilde{\bz} \circ \Phi$.
\par
$(2)$ If there exist diffeomorphism germs $\Phi:(I \times [0,2\pi),(t_0,\theta_0)) \to (\widetilde{I} \times [0,2\pi),(\widetilde{t_0},\theta_0))$ of the form $\Phi(t,\theta)=(\phi(t),\theta)$ and $\Psi:(\R^3,0) \to  (\R^3,0)$ of the form $\Psi(X,Y,Z)=(X,Y,\varphi(X^2+Y^2,Z))$ such that $\Psi \circ \bz=\widetilde{\bz} \circ \Phi$, then there exists a diffeomorphism germ $\psi: (\R^2,0) \to (\R^2,0)$ of the form $\psi(X,Z)=(X,\varphi(X^2,Z))$
such that $\psi \circ \gamma=\widetilde{\gamma} \circ \phi$
\end{theorem}
\demo
$(1)$ It is easy to see that $\Phi(t,\theta)=(\phi(t),\theta)$ and $\Psi(X,Y,Z)=(X,Y,\varphi(X^2+Y^2,Z))$ are diffeomorphism germs. 
By the assumption {of equivalence relation between curves}, 
we have {$(x(t),\varphi(x^2(t),z(t)))=(\widetilde{x}(\phi(t)),\widetilde{z}(\phi(t)))$}. 
Therefore, 
\begin{align*}
\Psi \circ \bz(t,\theta)&=\Psi(x(t)\cos \theta,x(t) \sin \theta,z(t))\\
&=(x(t) \cos \theta,x(t)\sin \theta, \varphi(x^2(t),z(t)))\\
&=(\widetilde{x}(\phi(t))\cos \theta, \widetilde{x}(\phi(t))\sin \theta, \widetilde{z}(\phi(t)))\\
&=\widetilde{\bz} \circ \Phi(t,\theta).
\end{align*}
$(2)$ By the assumption, $\phi$ and $\psi$ are diffeomorphism germs.
Since $\Psi \circ \bz(t,\theta)=\widetilde{\bz} \circ \Phi(t,\theta)$, we have 
$$
(x(t)\cos \theta, x(t)\sin \theta,\varphi(x^2(t),z(t)))=(\widetilde{x}(\phi(t))\cos \theta,\widetilde{x}(\phi(t))\sin \theta,\widetilde{z}(\phi(t))).
$$
It follows that $x(t)=\widetilde{x}(\phi(t))$ and $\varphi(x^2(t),z(t))=\widetilde{z}(\phi(t))$ and hence $\psi \circ \gamma(t)=\widetilde{\gamma} \circ \phi(t)$.
\enD

By the definition, we have the following result. 

\begin{proposition}\label{R-L}
Let $\gamma=(x,z), \widetilde{\gamma}=(\widetilde{x},\widetilde{z}):(I,t_0) \to (\R^2,0)$ be smooth curves.
If $x$ is $\mathcal{R}$-equivalent to $\widetilde{x}$ and $z$ is $\mathcal{L}$-equivalent to $\widetilde{z}$ at $t_0$, then $\gamma$ is equivalent to $\widetilde{\gamma}$ 
in the sense of Theorem \ref{diffeo-zero} $(1)$.
\end{proposition}

In detail for $\mathcal{L}$-equivalence see \cite{Bruce-Gaffney-dePlessis,Bruce-Giblin}.
If $\ord(f)(t_0)=k+1$ (see, (\ref{eq:order})) then $f$ is $\mathcal{R}$-equivalent to $\pm t^{k+1}$ (cf. \cite{Bruce-Giblin}).

\begin{proposition}
Let $\gamma=(x,z):(I,t_0) \to (\R^2,0)$ be a smooth curve.
If  {$x$ satisfies $\ord(x)(t_0)$ $=k+1$} and $z$ is regular at $t_0$, then $\gamma$ is equivalent to $t \mapsto (\pm t^{k+1},t)$ in the sense of Theorem \ref{diffeo-zero} $(1)$.
%
\end{proposition}
\demo
Since $x$ is $\mathcal{R}$-equivalent to $\pm t^{k+1}$ and $z$ is $\mathcal{L}$-equivalent to $t$, we have the result by Proposition \ref{R-L}.
\enD

\section{Examples}

We give concrete examples of surfaces of revolution with singular points.
\begin{example}{\rm
Let {$x: I\to\R$} be a function given by 
$$
x(t)=\exp\left(\int \cot{t} dt\right)=\exp\left(\int \frac{\cos t}{\sin t}  dt\right)=\exp (\log \sin t)=\sin{t},
$$ 
{where $I=(0,\pi)$. }
This is a solution of \eqref{differential-eq} when $\alpha(t)=-1$ and $\beta(t)=\cot{t}$. 
Thus we obtain the constant Gauss curvature $-1$ surface of revolution with singular point at $t=\pi/2$. 
By Theorem \ref{Gauss-curvature}, we have $z(t)$ as 
$z(t)=\cos{t}+\log\left(\tan ({t}/{2})\right)$, 
and hence the surface of revolution $\bz(t,\theta)=(x(t)\cos\theta,x(t)\sin\theta,z(t))$ 
of $\gamma(t)=(x(t),z(t))$ is a pseudo-sphere, see Figure \ref{fig:sing3} (cf. \cite{Gray}). 
In this case, we can take $\nu(t)=(\cos{t},-\sin{t})$ and $\bmu(t)=(\sin{t},\cos{t})$. 
Thus the curvature of $(\gamma,\nu)$ is $(\ell(t),\beta(t))=(-1,\cot{t})$. 
Since $\ord(\cos{t})(\pi/2)=\ord(\cot{t})(\pi/2)=1$, 
$\gamma$ actually has a $3/2$-cusp at $t=\pi/2$ (cf. Proposition \ref{constant-Gauss}). 
Moreover, since $\ell(t) \neq 0$ for all $t\in(0,\pi)$, 
one can consider at least an evolute $\mathcal{E}v^z(t,\theta)$ of $\bz(t,\theta)$ (cf. Proposition \ref{evolute} (2)). 
{On the other hand, 
$x(t)=\sin{t}\neq0$ on $I$ and $a(t)=\cos{t}$ vanishes at $t=\pi/2$. 
Thus the function $x(t)/a(t)=\tan{t}$ is defined on $\tilde{I}=I\setminus\{\pi/2\}$.}

We now take $\bn^z(t,\theta)$ as $(\cos{t}\cos{\theta},\cos{t}\sin{\theta},-\sin{t})$. 
Then an evolute $\mathcal{E}v^z(t,\theta)$ is given by
\begin{align*}
\mathcal{E}v^z(t,\theta)&=\bz(t,\theta)-\dfrac{\beta(t)}{\ell(t)}\bn^z(t,\theta)
=\left(\csc{t}\cos{\theta},\csc{t}\sin{\theta},\log\left(\tan\dfrac{t}{2}\right)\right).
\end{align*}
Setting $w(t)=\log\left(\tan({t}/{2})\right)$, then $dw/dt=\csc{t}\neq0$ for any $t\in(0,\pi)$. 
Thus this gives a parameter change and we have $t=2\arctan(\exp(w))$. 
Hence it follows that 
$$\mathcal{E}v^z(w,\theta)=(\cosh{w}\cos{\theta},\cosh{w}\sin{\theta},w).$$
This is a catenoid. 
{On the other hand, we can define another evolute $\widetilde{\mathcal{E}v^z}$ on $\tilde{I}\times[0,2\pi)$ by 
$$\widetilde{\mathcal{E}v^z}(t,\theta)=\bz(t,\theta)-\dfrac{x(t)}{a(t)}\bn^z(t,\theta)=\left(0,0,\log\left(\tan\dfrac{t}{2}\right)+\sec{t}\right).$$
Since $\tilde{I}\times[0,2\pi)$ is dense in $I\times [0,2\pi)$, the evolute $\widetilde{\mathcal{E}v^z}$ can be extended on $I\times[0,2\pi)$ continuously. 
This evolute $\widetilde{\mathcal{E}v^z}$ degenerates to the line (see Figure \ref{fig:sing3}).}

\begin{figure}[htbp]
  \begin{center}
    \begin{tabular}{c}

      \begin{minipage}{0.35\hsize}
        \begin{center}
          \includegraphics[clip, width=3cm]{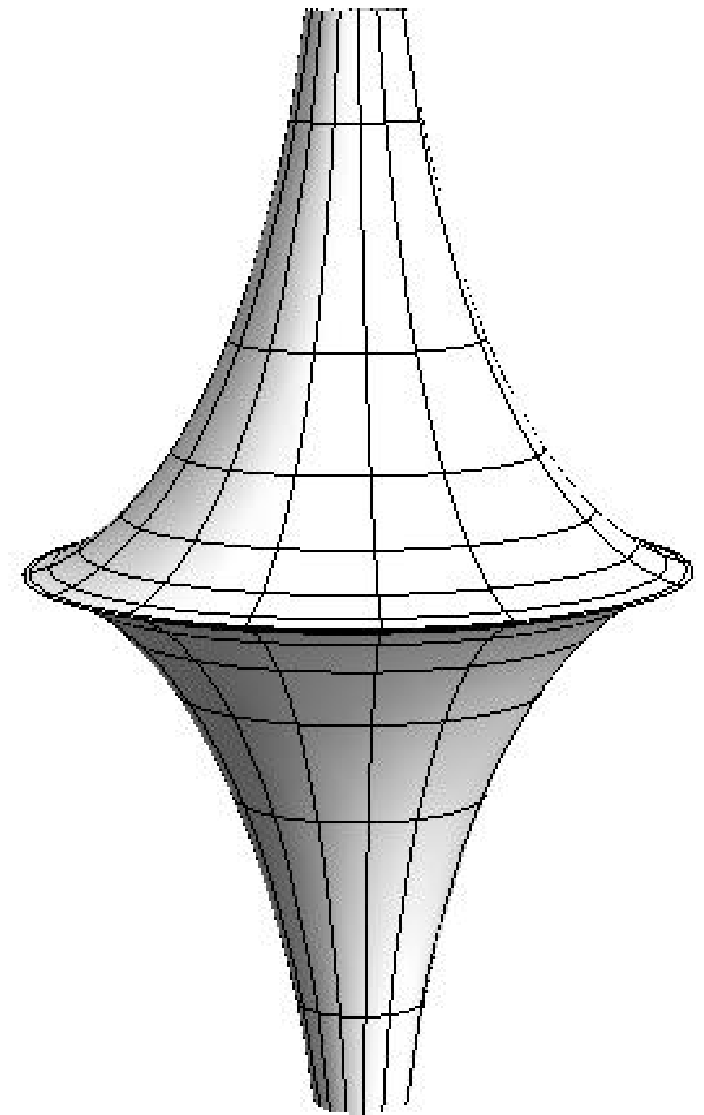}
        \end{center}
      \end{minipage}

      \begin{minipage}{0.35\hsize}
        \begin{center}
          \includegraphics[clip, width=3.5cm]{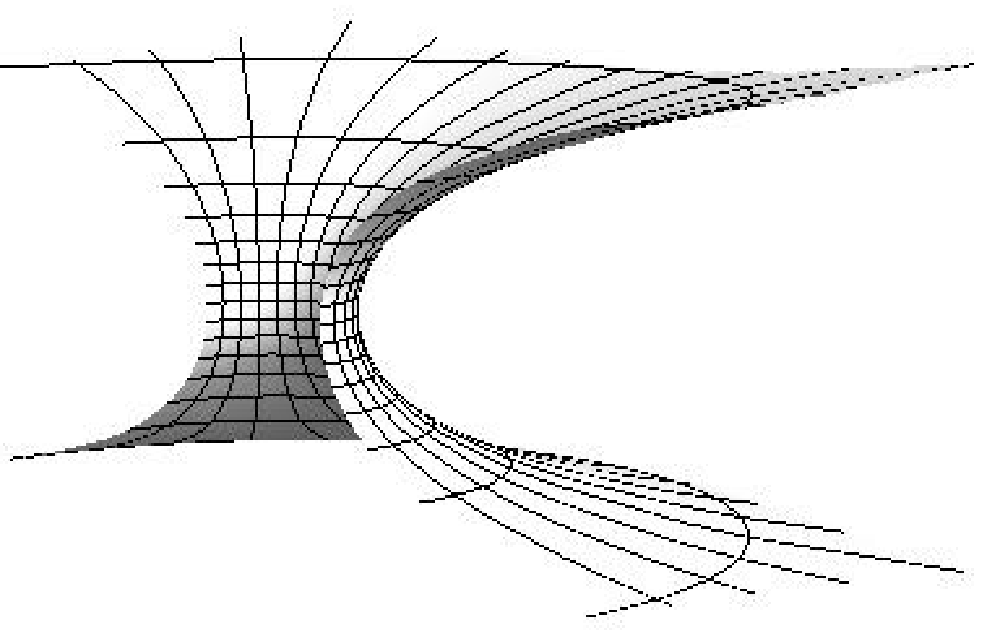}
        \end{center}
      \end{minipage}
\\
\begin{minipage}{0.35\hsize}
        \begin{center}
          \includegraphics[clip, width=4cm]{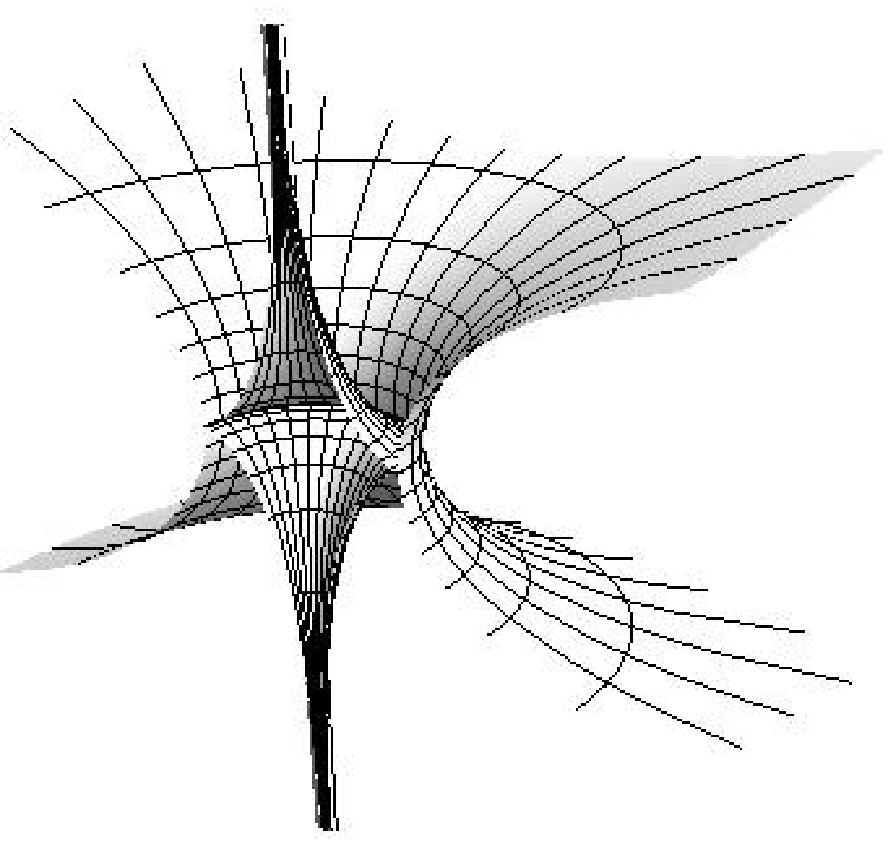}
        \end{center}
      \end{minipage}
  
\begin{minipage}{0.35\hsize}
        \begin{center}
          \includegraphics[clip, width=4.cm]{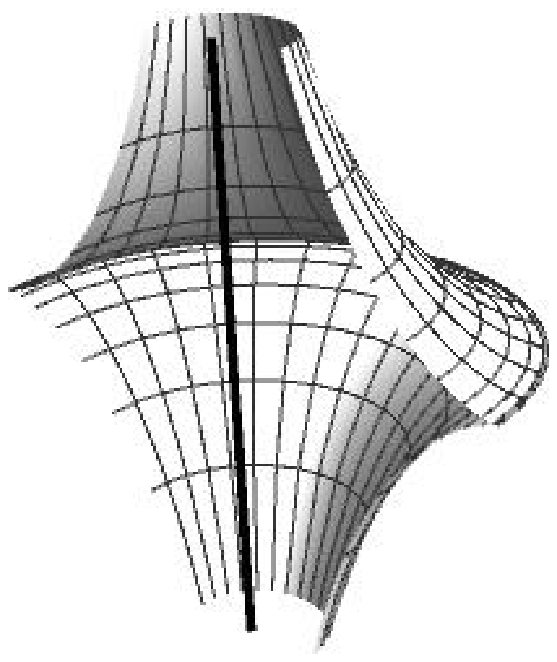}
        \end{center}
      \end{minipage}

    \end{tabular}
    \caption{From top-left to bottom-right: $\bz$, $\mathcal{E}v^z$, 
 both $\mathcal{E}v^z$ and $\bz$, both $\widetilde{\mathcal{E}v^z}$ (thick line) and $\bz$.}
    \label{fig:sing3}
  \end{center}
\end{figure}
}
\end{example}
\begin{example}{\rm
Let $\alpha(t)$ be a function satisfying $H_F^z=\alpha J_F^z$ for a surface of revolution $\bz(t,\theta)$ 
of the curve $\gamma(t)=(x(t),z(t))$. 
{Let $c_1$ and $c_2$ be initial values of $F(t)$ and $G(t)$, respectively.} 
Then $\gamma$ is given by \eqref{mean-xz}. 
We put $\alpha(t)=0$, $\beta(t)=t$ {($0\leq t<2\pi$)}, then $\gamma(t)$ and $\nu(t)$ are given by  
\begin{align*}
\gamma(t)&=\left(\sqrt{c_2^2 + \left(c_1 - {t^2}/{2}\right)^2},-c_2 \log\left(-2 c_1 + t^2 + \sqrt{4 c_2^2 + (-2 c_1 + t^2)^2}\right)\right),\\
\nu(t)&=\left(-\dfrac{c_1 -{t^2}/{2}}{\sqrt{c_2^2 + \left(c_1 - {t^2}/{2}\right)^2}},-\dfrac{c_2}{\sqrt{c_2^2 + \left(c_1 - {t^2}/{2}\right)^2}}\right)
\end{align*}
If we take $c_1=1/5$ and $c_2=3/10$, then we have a catenoid with singularity, see Figure \ref{fig:sing1}. 
\begin{figure}[htbp]
  \begin{center}
    \begin{tabular}{c}
      \begin{minipage}{0.33\hsize}
        \begin{center}
          \includegraphics[clip, width=4cm]{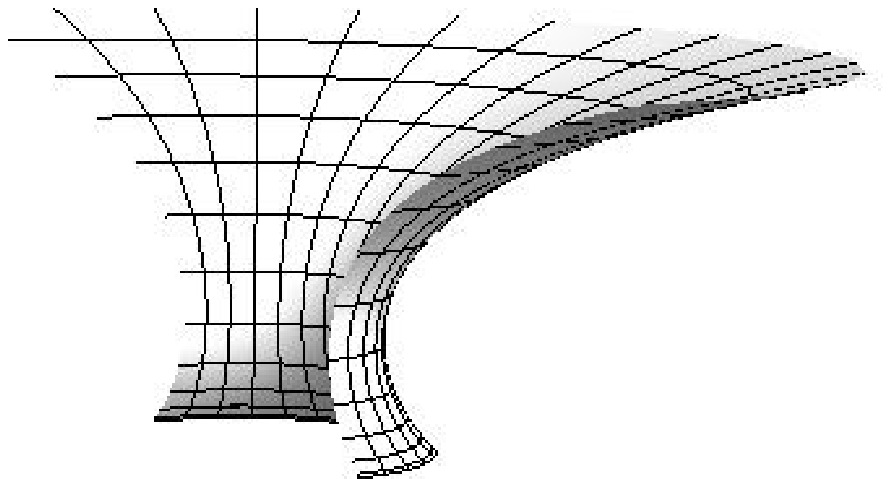}
        \end{center}
      \end{minipage}

      \begin{minipage}{0.33\hsize}
        \begin{center}
          \includegraphics[clip, width=4cm]{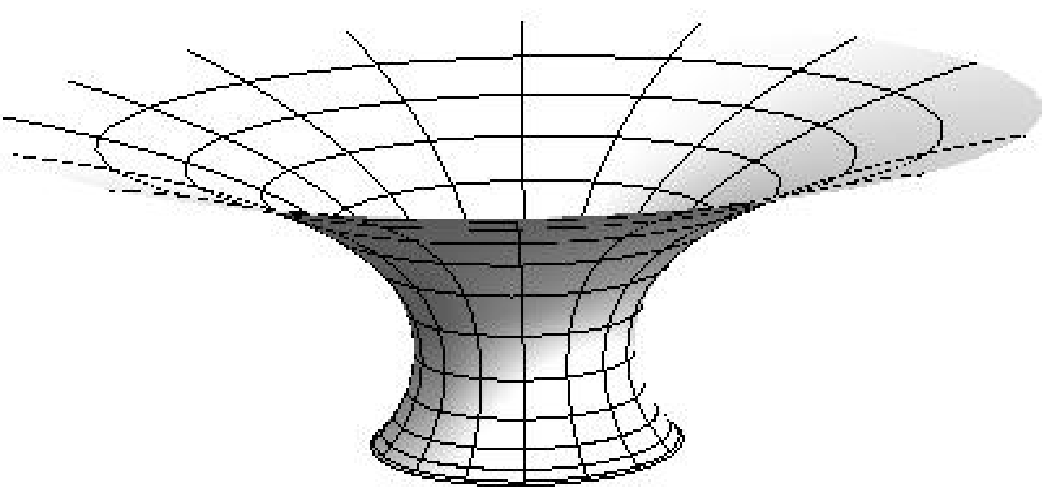}
        \end{center}
      \end{minipage}

    \end{tabular}
    \caption{Catenoid with singularities.}
    \label{fig:sing1}
  \end{center}
\end{figure}

We put $\alpha(t)=-1/2$ and $\beta(t)=t$ {($0\leq t <2\pi$)}. 
Then the surface of revolution is a constant mean curvature with singular point at $t=0$. 
In this case, $\gamma(t)$ is given by 
{
\begin{equation*}
\gamma(t)=\left(\sqrt{\left(c_1-\sin ({t^2}/{2})\right)^2+\left(c_2-\cos ({t^2}/{2})\right)^2}\right.,\\
\left.\int \dfrac{t}{x(t)}\left(1-2c_2\cos ({t^2}/{2})-2c_1\sin ({t^2}/{2})\right)dt\right),
\end{equation*}
}
{where $x(t)$ is the first component of $\gamma(t)$.}
We note that the integral of the $z$ component can be expressed concretely by using the Appell hypergeometric function $F_1$ (cf. \cite{Appell}). 
However it is too long, so we omit to write down here. 
Moreover, we see that 
{
\begin{equation*}
\nu(t)=\left(\dfrac{-1+c_1\sin ({t^2}/{2})+c_2\cos ({t^2}/{2})}
{\sqrt{\left(c_1-\sin ({t^2}/{2})\right)^2+\left(c_2-\cos ({t^2}/{2})\right)^2}},\right.\\
\left.\dfrac{-c_1\cos ({t^2}/{2})+c_2\sin ({t^2}/{2})}
{\sqrt{\left(c_1-\sin ({t^2}/{2})\right)^2+\left(c_2-\cos ({t^2}/{2})\right)^2}}
\right).
\end{equation*}
}
If we take $c_1=1/5$ and $c_2=3/10$, then we have an unduloidal surface (left-hand side of Figure \ref{fig:sing2}). 
On the other hand, if we take $c_1=1/5$ and $c_2=3/4$, we obtain a nodoidal surface (center of Figure \ref{fig:sing2}).
\par
Moreover, we put $\alpha(t)=t$, $\beta(t)=\sin{t}$ {($0<t<2\pi$)} and $c_1=c_2=1/10$. 
Then functions $F(t)$ and $G(t)$ are written as 
{
\begin{equation*}
F(t)=c_1 - \int \cos(2 (-t \cos t + \sin t)) \sin t dt, \quad
G(t)=c_2 - \int \sin(2 (-t \cos t + \sin t)) \sin t dt.
\end{equation*}
}
In this case, by $F(t)$, $G(t)$ above and $(1)$ of Proposition \ref{constant-mean}, we have 
the surface of  revolution with a $5/2$-cusp (at $t=\pi$) whose mean curvature is bounded, see right-hand side of Figure \ref{fig:sing2}. 

\begin{figure}[htbp]
  \begin{center}
    \begin{tabular}{c}

      \begin{minipage}{0.3\hsize}
        \begin{center}
          \includegraphics[clip, width=1.25cm]{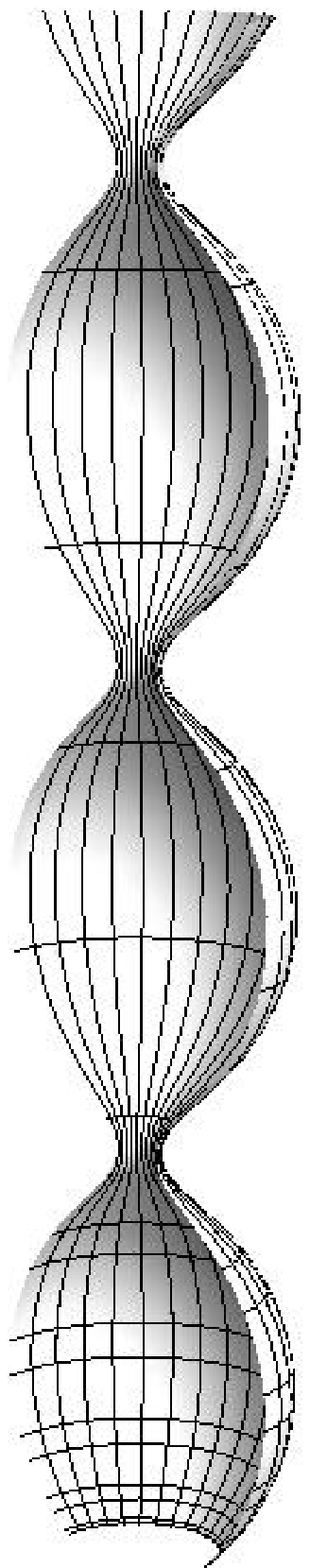}
        \end{center}
      \end{minipage}

      \begin{minipage}{0.3\hsize}
        \begin{center}
          \includegraphics[clip, width=3.5cm]{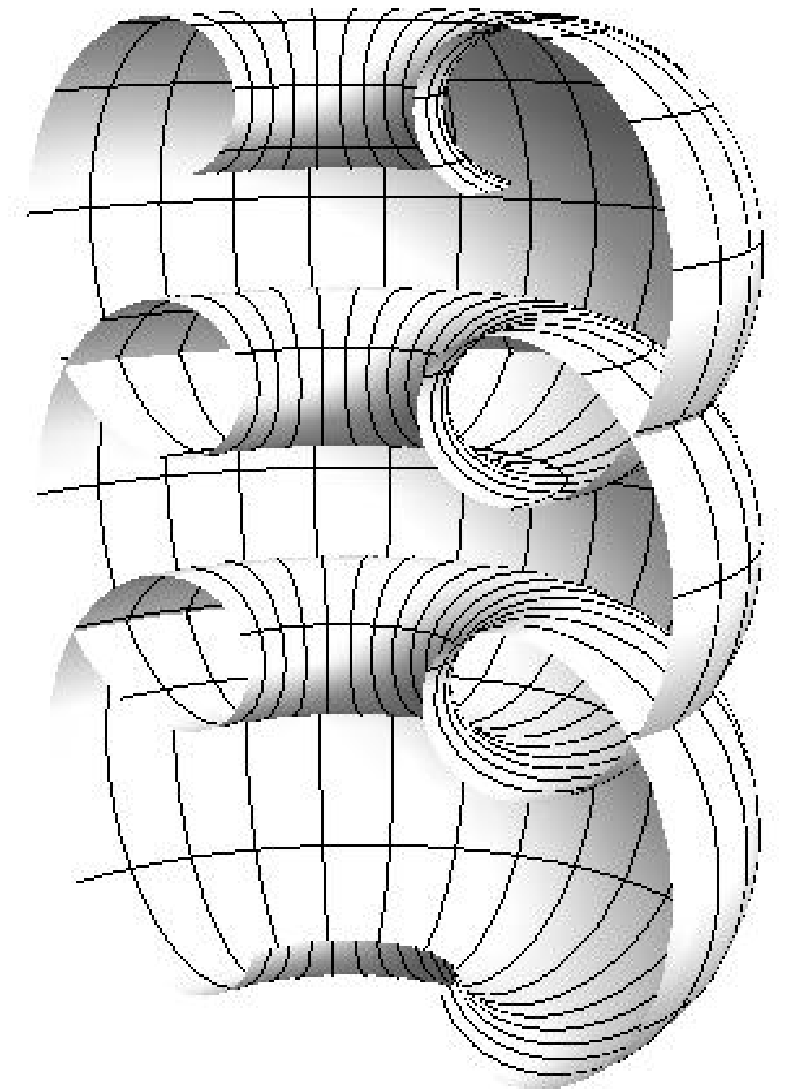}
        \end{center}
      \end{minipage}

\begin{minipage}{0.3\hsize}
        \begin{center}
          \includegraphics[clip, width=3.5cm]{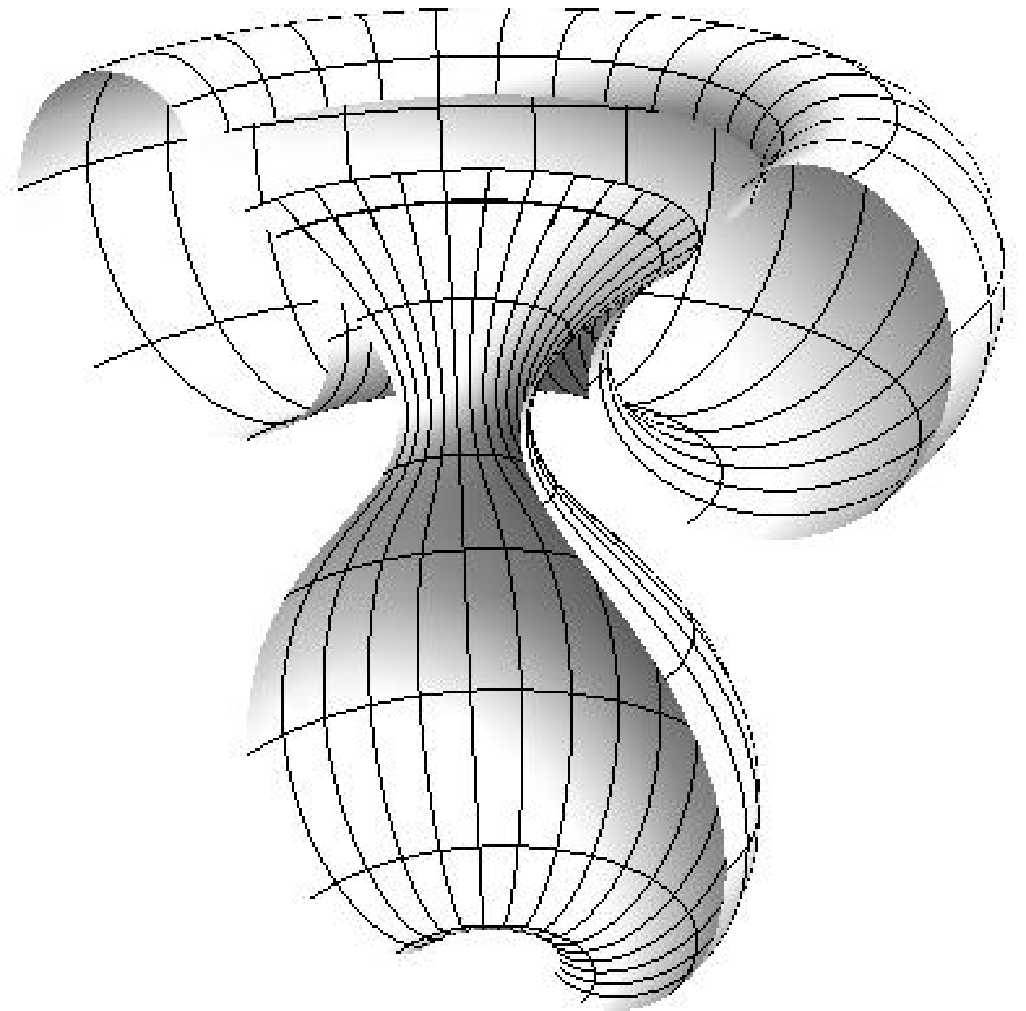}
        \end{center}
      \end{minipage}
    \end{tabular}
    \caption{Unduloidal surface (left), nodoidal surface (center) and surface of revolution with $5/2$-cusp (right).}
    \label{fig:sing2}
  \end{center}
\end{figure}
}
\end{example}

\ \\
Masatomo Takahashi, 
\\
Muroran Institute of Technology, Muroran 050-8585, Japan,
\\
E-mail address: \href{mailto:masatomo@mmm.muroran-it.ac.jp}{masatomo@mmm.muroran-it.ac.jp}
\\
\\
Keisuke Teramoto,
\\
Faculity of Science, Kobe University, Rokko 1-1, Kobe 657-8501, Japan
\\
E-mail address: \href{mailto:teramoto@math.kobe-u.ac.jp}{teramoto@math.kobe-u.ac.jp}

\end{document}